\documentclass[final,a4paper,10pt]{article}
\usepackage[brazil,english]{babel}     
\usepackage[latin1]{inputenc}
\usepackage{lmodern}
\usepackage[T1]{fontenc}
\usepackage{amscd,amsthm,amsmath,amssymb,latexsym}
\usepackage{mathrsfs,mathtools} 
\usepackage{indentfirst}

\newtheorem{teo}{Theorem}
\newtheorem{prop}{Proposition}
\newtheorem{lema}{Lemma}
\newtheorem{cor}{Corollary}
\newtheorem{ex}{Example}
\newtheorem{defini}{Definition}

\makeatletter
\newcommand{\subjclass}[2][1991]{%
  \let\@oldtitle\@title%
  \gdef\@title{\@oldtitle\footnotetext{#1 \emph{Mathematics subject classification.} #2}}%
}
\newcommand{\keywords}[1]{%
  \let\@@oldtitle\@title%
  \gdef\@title{\@@oldtitle\footnotetext{\emph{Key words and phrases: }#1.}: }%
}
\makeatother

\newcommand{\Addresses}{{
  \bigskip
  \footnotesize

  \medskip

   S.~O.~Juriaans,  
  \textit{E-mail address}: \texttt{ostanley@usp.br}

  \medskip

  J.~Oliveira, \textsc{Universidade Federal do Roraima 
Boa Vista  - RR - Brazil}.\par\nopagebreak
  \textit{E-mail address}: \texttt{joselito.oliveira@ufrr.br}

}}

\title{Fixed Point Theorems  for Hypersequences  and the Foundation of   Generalized  Differential Geometry I: The Simplified Algebra}

\author{ Juriaans, S.~O. \and Oliveira, J.}

\date{}

\begin{document}
	
\maketitle 
 
\begin{abstract}

Fixed point theorems are one of the many tools used to prove existence and uniqueness of  differential equations. When the data involved contains products of distributions,  some of these tools may not be useful. Thus rises the  necessity  to develop new environments and tools capable of handling such  situations.        J.F.  Colombeau,   E.E. Rosinger,     J. Jel\'\i nek,  M. Kunziger and J. Aragona,   together with their   collaborators,  were among the first to suggested, and propose,  a generalized differential geometry. The foundations of   a Generalized Differential  Geometry is set having Classical Differential  Geometry as a discontinuous sub-case,        a fixed point theorem for hypersequences is proved in the context of Colombeau Generalized Functions  and it is shown how it can be used to obtain  existence and uniqueness of differential equations whose data involve  products of distributions. Thus also setting the foundations of a Generalized Analysis.    The strain is also  picked  up  setting    the foundations  of  generalized manifolds   and   shown   that each classical manifold can be discretely embedded in a generalized manifold in such a way that the differential structure of the latter is a natural  extension of the differential structure of the former.  It is inferred    that $\ {\cal{D}}^{\prime}(\Omega)$ is discretely embedded in $\ {\cal{G}}(\Omega)$, that the elements of $\ C^{\infty}(\Omega)$ form a grid of equidistant  points in $\ {\cal{G}}(\Omega)$  and that association in $\ {\cal{G}}(\Omega)$ is a topological and not an algebraic notion. Ergo, classical solutions to differential equations are scarce.  These  achievements reckoned upon  the Generalized Differential Calculus invented by the first author  and his collaborators. Hopefully,     Generalized Differential Calculus and the   developments presented in this paper,  may  be of interest to those working in Analysis,  Applied Mathematics,  Geometry and  Physics.



\end{abstract}

\section{Introduction}

The Theory of Colombeau Generalized Functions is a nonlinear Theory of Generalized Functions which includes  Schwartz' Theory   of   Linear Generalized Functions, i.e., Schwartz Distribution Theory.    Colombeau's Theory   is well documented by now. Excellent textbooks and articles exist  that are pitstops  to appreciate and understand the theory and the wide spectra of applications. We  refer the reader to some of these excellent text,   \cite{JA1, AB, hebe, JFC1, JFC2, JFC3,  JFC4, grosser1, MK, MK1, MK2, MK3, NPS, Ni, OMM, OV, OV2, OV3,  OV4, OV5, PSV, scar,  EER1, EER2, EER3, EER4,  EER5, tod,  todver, V, V1},  to get the zest  of the basics  and relish    advanced parts   of this theory.

In the Colombeau environments, proving existence for differential equations involving products of distributions in  their data has always leaned  on classical results to guarantee existence. To achieve this, classical existence results are used,  proceeding  to prove  moderateness and  conclude existence of solutions in the environments of Colombeau Algebras. This can be highly nontrivial. One of the  setbacks  is that most tools used are not intrinsic  to these environments.  The development  of    Generalized Differential Calculus (see \cite{OJRO, OJR}) envisaged    the buildout of   tools,   intrinsic to the generalized environments,    making it possible to pin   less faith on the   classical ones. 

Let $M$ be an $n-$dimensional manifold.  The idea of linking a generalized objected $\widetilde{M}_c$  to  $M$ was first employed in \cite{MK4} where   a blueprint was given   how to use these objects to solve important problems in General Relativity.  Based  on this pivotal  idea, in \cite{OJR}, the notion of a generalized manifold was introduced.  The   definition  is exactly the same  as the classical one  the  difference being  that local charts take values   in open subsets of $\  \overline{\mathbb{R}}^n$ and differentiability is checked using the Generalized Differential Calculus.  In \cite{OJRO}, more details of Generalized Differential Calculus were worked out, showing that it extends and  behaves  very similar to classical Calculus and an example of a generalized manifold,   different from $\widetilde{M}_c$,  was  also given (actually it is  a subset of $\widetilde{M}_c$). As far as we know, other examples were not given yet and it remained unclear whether $\widetilde{M}_c$ was a generalized manifold and whether there existed other examples. 

A pursue in another direction  was the construction of an diffeomorphism  invariant Colombeau algebra. This was early  undertaken in \cite{colmer, jel} and was settled in the definite in \cite{MK3} well afore  Generalized Differential Calculus was proposed.   These are top-notch  papers which  show that  the main  obstruction to the  construction of   such an algebra has a   topological nature:   Colombeau algebras  are ultrametric  spaces which naturally mismatch with the classically used topologies. This translate into  a highly  non-trivial endeavor   the creation of  an algebra that can be attached to classical manifolds. The last example given in \cite{MK3} shows that  having such an algebra does not necessarily make  things much easier  when applying the theory to obtain existence of solutions of  differential equations having products of distributions in their data. It is essential to observe  that  an algebra of generalized functions that can   be attached to a manifold was also achieved in \cite{rossing}. Amazingly enough, in this case, technicalities are not that involved. 

  In \cite{OJR}, all necessary machinery of Generalized Differential Calculus (such as the Inverse Mapping Theorem, The Implicit Function Theorem and others) were proved so that a consistent basis could be laid  for a Generalized Differential Geometry.  At first,    definitions  given  and results obtained  are exactly the same  as the classical ones   but  extend   the latter  in a non-trivial way.  Howbeit,   much has yet to be accomplished before   this Generalized Differential  Geometry  unveils   its smoldering potential.   Generalized  Differential  Calculus   is based on key ideas developed over the years by all prominent researchers in the field but the decisive ideas are   due to Kuzinger-Oberguggenberger (\cite{MK1}) and Biagioni-Scarpal\'ezos (\cite{hebe, scar}).  The topology in use (see \cite{top1, top2, AGJ})  is a slide modification of the sharp topology   introduced by Biagioni-Scarpel\'ezos (see \cite{hebe, scar}), yet equivalent to it,  is more natural and in harmony with the algebraic structure (see for example \cite{AGJ, JRJ}) of the Colombeau algebras.   An interesting fact   is that,   in the sharp  topology,  $\mathbb{R}^n$  embeds as  a discrete subset of $\ \overline{\mathbb{R}}^n$ and yet the Generalized Differential Calculus is a near perfect  extension of the Newtonian Differential Calculus (see the Embedding Theorem in \cite{OJR}). In particular,    Classical Space-Time becomes a grid of equidistant  points  in   Generalized Space-Time,  a possibility that was   raised along time by many physicists and  more recently also  in \cite{wolfram}.  The common   distance betwixt    grid points  could  well be glossed as Planck's constant.  
  
  Can   this discontinuity in classical space-time  be perceived experimentally? Or at least, can one be convinced that we do have an issue  in this direction?   This is where the notion of  hypersequences steps into the picture. Since classical  space-time is a grid of equidistant  points  it is impossible for  classical sequences  to  converge  in this new environment. In particular,  it is no longer true that the sequence $(\frac{1}{n})_{n\in \mathbb{N}}$  converges  in the ring of  Colombeau  generalized numbers. But the hypersequence it generates is of the form $(\frac{1}{n})_{n\in\widetilde{\mathbb{N}}}$  and does  converge to zero in the Colombeau environment.   From  the point of view of someone living in generalized space-time, classical convergence of a sequence $(x_n)$  is equivalent to the existence of a $n_0\in \widetilde{\mathbb{N}}$ such that $x_n -x_m\in V_1(0)$ if $n,m >n_0$.     So classically we only measure upto scale $\alpha=[\varepsilon \longrightarrow \varepsilon]$,  which  is the reason  for  calling it  our {\it natural gauge},  the latter being  first introduced  in \cite{JO}.  A similar problem occurs when   proving  existence of differential equations    using  classical tools to prove moderateness  and existence in    environments  that are like chalk and cheese, in the topological sense, compared to the environment where these tools come from. It takes  much more to be convincing.

The paper  is  structured as follows.  In the next section we recall the necessary machinery needed to understand the context and prove subsequent results. In the third section we prove a fixed point theorem for hypersequences, prove  that association is  a topological and not an algebraic concept  and that ${\cal{D}}^{\prime}(\Omega)$ is discretely embedded in ${\cal{G}}(\Omega)$ thus proving that classical functions are extremely rare.  In the fourth   section we prove that $\widetilde{M}_c$ is a generalized manifold  and devote  the last section to examples and the enumeration of  some  results  in this new  generalized geometry.  This is the first of two papers.  The second paper is in the context of the full algebra of Colombeau Generalized Functions  (see \cite{paulo}) thus completing the  proposal of this new Generalized Differential  Geometry as a roundabout route to define generalized functions on  manifolds.   

The notation $\overline{\mathbb{K}}$, for the  ring of Colombeau generalized numbers,  was introduced by Colombeau. However,   developments overtime show that it is   more reasonable to use the notation  $\widetilde{\mathbb{K}}$  to denote this ring. Here we will still be using Colombeau's original notation to be consistent with the notation in \cite{OJRO, top1, top2, OJR}.  

This paper was written while the second author held a   pos-doc position at IME-USP, the University of S\~ao Paulo-Brazil. The dimension invariance  theorem of section four  and some of the examples   of the  last  section  are part of his Ph.D.        thesis (\cite{jose})   written under the supervision  of the first author.

\section{Preliminaries}

We shall mainly  work over the field $\ \mathbb{R}$ of real numbers but all results also hold  for $\mathbb{C}$.  This is the reason why sometimes we use $\mathbb{K}$ to   denote either of these fields.  One could  rightfully ask "why  not consider the field $\ \mathbb{Q}\ $?"  The answer is simples:  The Colombeau Theory constructed using $\ \mathbb{R}$ is the same as the Colombeau Theory constructed using $\ \mathbb{Q}$ since real numbers can be seen as nets of rational numbers. The reason why we end up with a bigger structure,  which is not a field but  never the less  very interesting, is because we  mod out some, but not all,  nets converging to zero. These surviving nets,  converging to zero,  are the infinitesimals which inhabit the halos of the elements of the newly formed environments.    

Set $I=]0,1]$, $I_{\eta}=]0,\eta]$, for $\eta\in ]0,1[$  and let   $\alpha$  be the identity map $\alpha : I\longrightarrow \mathbb{R}$, $\alpha(\varepsilon)=\epsilon$.  We shall denote, once in a while, $\alpha_n =\alpha^n$ and call $\ \alpha$  the {\it standard} or {\it natural gauge}. Nearly all results in this paper can be proved for other gauges using the already existing results for these gauges (see \cite{rob, V2}). 

  A map (also called a net)  $x : I\longrightarrow \mathbb{R}$ is moderate if $|x| <\alpha^r$, for some $r\in \mathbb{R}$, i.e. $|x(\varepsilon)| <\varepsilon^r$, $\forall \varepsilon \in I_{\eta}=\ ]0,\eta], \exists \ \eta <1$.  Denote the set of moderate maps by ${\cal{E}}_M$ and  by ${\cal{I}}=\{ x : x \ \mbox{is moderate and }\ |x|<\alpha^n, \forall n \in \mathbb{N}\}$. For  $x\in {\cal{E}}_M$,  denote by $V(x)=Sup\{r\in \mathbb{R} : |x|<\alpha^r\}$ and set $\|x\| =e^{-V(x)}$.  Then $\cal{I}$ is a radical  ideal of the ring  ${\cal{E}}_M$ and setting, $\overline{\mathbb{R}} : =\frac{{\cal{E}}_M}{\cal{I}}$, we have that $(\overline{\mathbb{R}} , \| \ \|)$ is an ultrametric partially ordered topological  ring whose group  of units, $Inv(\overline{\mathbb{R}} )$, is open and dense (see \cite{JO}). The latter property is essential in developing the Generalized Differential Calculus (\cite{OJRO, OJR}).  This  topological ring, $(\overline{\mathbb{R}} , \| \ \|)$, is called the ring of Colombeau Generalized (real) numbers.  A generalized number  $x\in \overline{\mathbb{R}} $ is a unit if and only if $|x| >\alpha^r$ for some  $r\in \mathbb{R}$   and it is a non-unit if and only if there exist a nontrivial  idempotent $e\in {\cal{B}}(\overline{\mathbb{K}})\ $ such that $e\cdot x=0$ (see \cite{JO, JJO}). In particular, a generalized number is either a unit or a zero divisor. The ring $\overline{\mathbb{R}}$, contains  $\ \mathbb{R}$ as a discrete subfield. Actually,  $\ \mathbb{R}$  is a grid of equidistant  points in $\overline{\mathbb{R}}$.   The latter  is a partially ordered ring whose maximal ideals and idempotents have been completely determined (see \cite{OJR, JO, JJO, V}). The partial order is not intrinsic  but stems from the order of $\ \mathbb{R}$.  This is maybe the only definition that is not, yet,  intrinsic. Distance emerges  from this order and that is why it is important to understand order.   In the   references we just mentioned,  one finds the following facts:  the Jacobson Radical of $\ \overline{\mathbb{R}}$ is trivial, its ideals are convex, its Krull dimension is infinite and it has a minimal prime  which is also a maximal ideal. Its Boolean algebra, ${\cal{B}}(\overline{\mathbb{R}})$,  consists of $\{0,1\}$ and  positive elements each  of which is a characteristic functions of a subset $S\subset I$, such that $0\in \overline{S}\cap\overline{S^c}$, where the last two  bars stand for  topological closure in $\ \mathbb{R}$. The set of these subsets is denoted by $\ {\cal{S}}\ $  and was defined in \cite{JO}. Ultrafilters of $\ {\cal{S}}$ partially parametrize prime and maximal ideals of $\ \overline{\mathbb{K}}$.   It also holds that ${\cal{B}}(\overline{\mathbb{R}})={\cal{B}}(\overline{\mathbb{C}})$ (see \cite{JJO}). In particular, the Heaviside function $H\notin {\cal{B}}(\overline{\mathbb{R}})$, i.e., $H^2\neq H$ (see \cite{jailson}).   

H. Biagioni and Scarpel\'ezos were the first to suggest  the  topology,  defined above,  for $\overline{\mathbb{R}}$. It   came to be known as the {\it sharp topology} turning  $\overline{\mathbb{R}}$ into a complete ultrametric algebra and hence, its topology is generated by balls. In \cite{top1, top2} it was shown  that this topology was also generated by the sets $V_r[x]=\{y\in \overline{\mathbb{R}}: |y-x|<\alpha^r\}$, balls, with generalized numbers as radii ,   compatible with the ring structure.  It is easily seen  that $B_{2r}(0)\subset V_r[0]\subset B_{r/2}(0)$ if $r>0$.
 
Let  $\ \Omega \subset \mathbb{R}^n$ be an open subset    with  an exhaustion by relatively compact  subset  $\Omega_m\subset \Omega$.  Consider   nets $p=(p_{\varepsilon}), p_{\varepsilon}\in  \Omega_{m(p)}\ \forall \varepsilon\in I_{\eta}$, $m(p)\in \mathbb{N}$,   such that  the net $(\|p_{\varepsilon}\|)\in {\cal{E}}_M$.  Factoring out   nets  $p$ for which also   $(\|p_{\varepsilon}\|)\in \cal{I}$,   give rise to  a subset of  $\ \overline{\mathbb{R}}^n $ which is  denote by $\widetilde{\Omega}_c$ (see \cite{MK, MK1}).  The notation $\widetilde{\Omega}$ is used if  one does not require the existence of $m(p)$.  The algebra of Colombeau generalized functions ${\cal{G}}(\Omega)$ (see \cite{AB, JFC1, JFC2, MK} for the original definition), defined on the open subset $\Omega\subset \mathbb{R}^n$,  can be viewed as $C^{\infty}-$functions defined  on $\widetilde{\Omega}_c\subset \overline{\mathbb{R}}^n$ and taking values in $\ \overline{\mathbb{R}}$ (see \cite{MK, MK1, OJRO, OJR}).  In \cite{OJR} (see also \cite{OJRO}) the foundation of  Colombeau Generalized Calculus is laid and shown that  ${\cal{G}}(\Omega)$  can be embedded into $C^{\infty}(\widetilde{\Omega}_c,\overline{\mathbb{R}})$. In particular, Schwartz space of linear distributions,   $D^{\prime}(\Omega)$,   can be seen as infinitely differentiable functions  where differentiability is defined a la Newton. So we have come full circle from seeing  elements of  $\ D^{\prime}(\Omega)$ as linear maps, and hence not undergoing variation,  to seeing them as functions undergoing variation (note that, classically, derivation in $\ D^{\prime}(\Omega)$ is defined without the use of variation).  (An  interesting fact  is that in the proof of the differentiability of objects in  the construction  of an invariant algebra in \cite{MK3},  Schwartz distributions are treated as linear maps. In the proof of  their embedding in this newly constructed algebra the linearity of a Schwartz distribution  is not used. When a Differential Calculus for these invariant algebras is ready, distribution will have both the classical and generalized meaning). More on this, in the next section (see also the remark after \cite[Proposition 1.6.3, 1.7.28]{MK} and the remark after the latter). For example, the Delta Dirac function $\delta$ is the derivate of $\ H$, the Heaviside function, and the calculation is just the ordinary calculation from classical Calculus. Another interesting point is that,  in the presence of  moderateness,  negligibility only has to be checked at level $0$.  This is mentioned and proved in several references. See, for example, \cite[paragraphs after I.Theorem 7.13]{MK3}. Generalized Differential Calculus allows to give an easy proof of this fact.  In fact, in the presence of moderateness negligibility at level $0$  is a statement about point values: If   $\hat{f}(\varepsilon, x) $ is moderate, then it defines an element   $f\in \ C^{\infty}(\widetilde{\Omega}_c,\overline{\mathbb{R}})$. Given $x\in \widetilde{\Omega}_c$ there exists a compact subset $K\subset \Omega$ containing a representative of $\ x$.   Moderateness at level $0$ implies that $ \|(f)_{|_K}\|_{\infty}=0$ (note that this is exactly  the uniformity on $\ K$). It  follows that  $f(x)=0$, proving that   $\ f=0$ in $\widetilde{\Omega}_c$. Generalized Differential Calculus gives that $\partial^{\beta}f=0, \forall\ \beta \in \mathbb{N}^n$. Since the Embedding Theorem \cite[Theorem 4.1]{OJR} tells us that derivations commutes with the embedding,  it follows that $f=0$ in ${\cal{G}}(\Omega)$.   So there is no need to check other levels. The same proof holds for the full algebra and should also work for the invariant algebra once we have at hand  a Generalized Differential Calculus for the latter. Note however that one must have negligibility at level $0$ and not just point values of elements of $\ \Omega$  being zero. To see why, consider $f=x\delta\in{\cal{G}}(\mathbb{R})$, where $\delta=[(\rho_{\varepsilon})]$, $\rho$ a mollifier.  We have that $f(x)=0, \forall x\in \mathbb{R}$.  Also, for $\ x_0\in \mathbb{R}$,   $f(x_0\alpha)=x_0\cdot\rho(x_0)$,   for  $\varphi\in {\cal{D}}(\Omega)$ we have that   $\int\limits_{\mathbb{R}}f(x)\varphi(x)dx=[(\int\limits_{\mathbb{R}}(x\varphi)(x)\rho_{\varepsilon})]$ and hence, since $(\rho_{\varepsilon})$ is a delta-net,  
$\int\limits_{\mathbb{R}}f(x)\varphi(x)dx  = (\int\limits_{\mathbb{R}}x\varphi(x)\rho_{\varepsilon}(x)dx )    \approx (x\varphi)(0)=0$. This proves that $\  f\approx 0$ but $f\neq 0$.  This example also shows an interesting phenomena: $f(0)=0$ and  thus  for $x\in V_r(0)$,  a small enough sharp neighborhood of $\ 0$, $f(x)\in V_1(0)$. For classical mathematics (and hence measurements) $f(x)=0$ and thus seemingly does not interfere with physical reality. But for histories of the form  $x=x_0\alpha, x_0\in \mathbb{R}$ we have that $f(x)=x_0\cdot\rho(x_0)\in \mathbb{R}$ and thus interfere with physical reality. These "waves" of appearing and disappearing from physical reality  are the source of the turbulence effects we see in physical reality. And it can be worse. Consider  $g(x)=f^k$,  $k\in \mathbb{N}$. Then $g(0)=0, g(\alpha)=(\rho(1))^k$. So if $\rho(1)>1$ and $k$ is large then these "waves" coming from $V_r(0)$  which   we cannot measure,  can effect in  a non-trivial way physical reality.  Note also that $\delta(0)=\alpha^{-1}\cdot \rho(0)$ is an infinity we can not measure but it is cancelled out on histories, $x_0\alpha$  near $0$,  giving us a real number,  $f(x_0\alpha)=x_0\cdot\rho(x_0)\in \mathbb{R}$,  that we can measure. Even though the history $x_0\alpha$  is near $0$, the position in physical reality where we observe the effect can be faraway from $0$ (in this case at $x_0\in \mathbb{R}$) and the result of the measurement $x_0\rho(x_0)$ becomes small as $x_0$ goes to infinity.     So turbulence should  be the interaction of elements of $\ B_1(0)$ and infinities, i.e.,  elements of norm greater then $1$, producing a measurable but not predictable effect on physical reality. The nonpredictability stems  from the fact that spheres in Colombeau environments are clopen  sets  and classical space-time is a grid of equidistant points. Jumps from one sphere to another sphere  occur multiplying with the $\ \alpha^r$, $\ r\in \mathbb{R}$, which form  a discrete  chain of  quanta.

The construction carried out above with $\Omega\subset \mathbb{R}^n$ can also be carried out with any subset $X\subset \mathbb{R}^n$. In fact,  consider $\ X$ with the induced topology,  consider an exhaustion $(X_n)$ of $\ X$ by relatively compact subsets and proceed as before. The set obtained in $\overline{\mathbb{R}}^n$ will be denoted by $\widetilde{X}_c$. We embed $X$ into $\ \widetilde{X}_c$ using constant nets $\ p=(p_{\varepsilon}), p_{\varepsilon}=x\in X,  \varepsilon \in I_{\eta}$. It is clear that, in the sharp topology, $X$ is a discrete grid of equidistant  points contained in $\widetilde{X}_c$. The remarkable thing is that, in case of a submanifold  $ \ M$  of $\  \mathbb{R}^n$, Generalized Differential Calculus on $\widetilde{M}_c$ will be a generalization of the Classical Differential Calculus on $\ M$, although $\ M$ is discretely embedded in $\ \widetilde{M}_c$.

Let $q$ be any norm on $\mathbb{R}^n$. Extending it in the obvious way to $\overline{\mathbb{R}}^n$, we define for     $ x=(x_1,\cdots, x_n)\in \overline{\mathbb{R}}^n$,  $\|x\|_q =q(x)\in \overline{\mathbb{R}}$ and  ${}_q\|x\|\in \mathbb{R}$ to be the norm of $q(x)$ as an element of $\ \overline{\mathbb{R}}$.  If $q(x)=\sqrt{x_1^2+\cdots +x_n^2}$ then we  write $\|x\|_q=\|x\|_2$ (see \cite{OJRO, JJO}).  Since all  norms on $\mathbb{R}^n$ are equivalent, it is easily seen that ${}_q\|x\|$ does not depend on the norm $q$ and thus we shall write it as $\|x\|$.

The positive cone, $\  \overline{\mathbb{R}}_+$,  of  $\  \overline{\mathbb{R}}$ is not an open subset. In  fact, let $e$ be an  idempotent and set $\ x_n =e-(1-e)\cdot \alpha^n$. Then $|x_n|= e+(1-e)\cdot \alpha^n$.  We clearly have that $x_n$ is not in the positive cone but  $x_n\longrightarrow e$.   However if we let  $Inv( \overline{\mathbb{R}})_+=Inv( \overline{\mathbb{R}})\cap \overline{\mathbb{R}}_+$  then we have:

\begin{lema}  Let  $Inv( \overline{\mathbb{R}})_+=Inv( \overline{\mathbb{R}})\cap \overline{\mathbb{R}}_+$. 

\begin{enumerate}

\item  $Inv( \overline{\mathbb{R}})_+$ is  an  open subgroup  of $\  \overline{\mathbb{R}}$.
 
 \item Let $t\in \widetilde{[0,1]}$ and $x,y\in Inv( \overline{\mathbb{R}})_+$. Then $tx+(1-t)y\in Inv( \overline{\mathbb{R}})_+$.
 \end{enumerate} 
 
\end{lema}
 
 \begin{proof} The fact that   $Inv( \overline{\mathbb{R}})_+$ is a subgroup of $\  \overline{\mathbb{R}}$ is clear. So 
 let $x\in Inv( \overline{\mathbb{R}})\cap \overline{\mathbb{R}}_+$.  Since $x$ is invertible,  there exists $\alpha^r$ such that $x>\alpha^r$. If  $y\in V_r(x)$ then $|y-x|<\alpha^r$ and thus $0<x-\alpha^r <y$. On the other hand, since $Inv( \overline{\mathbb{R}})$ is open and the $\alpha_t$'s form a totally ordered set, we may take $r$  such that $V_r(x)\subset Inv( \overline{\mathbb{R}})$.
 
 To prove the second part, take $m>0$ such that $\alpha^m <min\{x,y\}$. Then if follows that $tx+(1-t)y\geq t\alpha^m +(1-t)\alpha^m=\alpha^m$ and thus, $tx+(1-t)y\in Inv(\overline{\mathbb{R}})$.
 \end{proof}

The negative cone $Inv( \overline{\mathbb{R}})_+=Inv( \overline{\mathbb{R}})\cap \overline{\mathbb{R}}_{-}$ is also an open subset of $\  \overline{\mathbb{R}}$. It follows that $Inv( \overline{\mathbb{R}})$ has two connected component which are both open subsets of $\ \overline{\mathbb{R}}$. Moreover,  $Inv( \overline{\mathbb{R}})_+\cap Inv( \overline{\mathbb{R}})_- =\emptyset$  and $0$ is in   the topological closure of both $Inv( \overline{\mathbb{R}})_+$ and $Inv( \overline{\mathbb{R}})_-$. Both are  closed under addition and interleaven (see the end of this section or \cite{OV}).

  The lemma shows that there can not exist  a continuous curve whose initial value is negative, its final value is positive and at all other instants  its values are  comparable with $0$.  This is   exactly what is needed to define the notion of orientation on generalized manifolds.

The  ideas to consider nets of point in $\mathbb{R}^n$ was  introduced in \cite{grosser1, MK, MK1}. This was used  in \cite{OJRO} to define  the notion of  membranes and histories in $\overline{\mathbb{R}}^n$.   Subsequently,  in \cite{OV},  the notions  of internal and strong internal sets (internal sets are generalization of  membranes)   were  introduced, inspired also by concepts    of  nonstandard analysis.  In this same paper,  (\cite{OV}), very strong and relevant   properties  involving these notions were proved.    For example, it is proved that strongly internal sets are  open subsets of $\ \overline{\mathbb{R}}^n$ where as  internal sets are closed subsets  of the same space.  We will be using freely the results contained in these references.

Given  a net $(A_{\varepsilon})$ in $\mathbb{R}^n$ we shall write it  also as $A_{\alpha}$, being $\alpha=[\varepsilon \rightarrow \varepsilon]$ our natural or standard gauge.  For an idempotent $e\in {\cal{B}}(\overline{\mathbb{K}})$ we define $e\alpha=e(\varepsilon)\varepsilon$, meaning that when $e(\varepsilon)=0$ this index will be omitted. We also write $eA_{\alpha}$ for the net $A_{e\alpha}$,  $\partial A_{\alpha}$ for the net $(\partial A_{\varepsilon})$ (the boundary)   and     $int(U)$ for the set of interior point of a subset $\ U\subset \mathbb{R}^n$. Given a net $A_{\alpha}\subset \mathbb{R}^n$ of subset of $\ \mathbb{R}^n$, we denote the membrane, or internal set,  it originates by $[A_{\alpha}]$ and the strongly internal set it originates by $\langle A_{\alpha}\rangle$. As mentioned above,   internal sets are closed in the sharp topology while strongly internal sets are open in the sharp  topology.  We say that $(A_{\alpha})$ is {\it regular} if there exists $k\in \mathbb{N}$ such that at each boundary point $A_{\varepsilon}$ we can inscribe  balls of radius $\varepsilon^k$ tangent to $ \partial A_{\varepsilon}$ and contained in $\ int(A_{\varepsilon})$ and another  ball of the same radius tangent at  the same point but contained in $\ ( int(A_{\varepsilon}))^c$.   As a result,   the volume of a regular net is a unit since $vol(A_{\alpha })\geq vol(V_k(0))=\pi\alpha^{2k}\in Inv(\overline{\mathbb{K}})$ (see \cite{OJRO}). For example, this is the case if the boundaries are compact   hyper surfaces whose encompassing volumes are not shrinking too fast.

\begin{lema}
Let    $(A_{\alpha})$ be a  net in  $\mathbb{R}^n$  and $U=\langle A_{\alpha}\rangle$ its strong internal set. Then  $\partial \langle A_{\alpha}\rangle =[ \partial A_{\alpha}]$. 
\end{lema}

\begin{proof}
For  $z\in \partial \langle A_{\alpha}\rangle$,  there exist sequences $(z_n)\subset \langle A_{\alpha}\rangle$ and $(p_n)\subset  (\langle A_{\alpha}\rangle)^c$ both converging to $z$ in $\overline{\mathbb{K}}^n$. Let $w=dist(z,[\partial A_{\alpha}])$ be the distance of $\ z$ to the membrane $[\partial A_{\alpha}]$. If $w=0$  then $z\in [ \partial A_{\alpha}]$. If not, then there exist $e\in {\cal{B}}(\overline{\mathbb{K}})$ and $t > >k$ such that $e\cdot w >e\cdot \alpha^t$. In particular, $dist(z_{\varepsilon},\partial A_{\varepsilon})>\varepsilon^t \neq 0$ if $\ e(\varepsilon)=1$.  Hence, for $\varepsilon$ such that $e(\varepsilon)=1$, we have that either $z_{\varepsilon}\in  int(A_{\varepsilon})$ or $z_{\varepsilon}\in int((A_{\varepsilon})^c )$. Consequently, we may write $e=e_1+e_2$,  a sum of orthogonal idempotents,  such that $e_1\cdot z\in e_1 (int( A_{\alpha}))$ and $e_2\cdot z\in  e_2(int(A_{\alpha})^c)$.  On the other hand, since $p_n\rightarrow z$ , there exists $n_0$ such that if $n>n_0$ we have that $dist(e_1\cdot p_n,e_1\cdot [(\partial A_{\alpha})]) >e_1\cdot \alpha^{7t}$. But since $e_1\cdot z\in e_1\cdot\partial\langle A_{\alpha}\rangle$ this implies that $e_1\cdot p_n\in \langle e_1\cdot A_{\alpha}\rangle$, a contradiction, unless $e_1=0$. If so, then revers the roles of $\ z_n$ and $p_n$ obtaining another contradiction.  Thus we have that $w=0$ and the result is proved. \end{proof}


In case $A_{\alpha}$ consists of intervals $J_{\varepsilon}=]a_{\varepsilon},b_{\varepsilon} [\subset  \mathbb{R}$, uniformly bounded,  we have that $\partial \langle A_{\alpha}\rangle= Interleaven\{a,b\}$, where $a=[(a_{\varepsilon})] $ and $b=[(b_{\varepsilon})] $.   The notion of interleaven is defined in \cite{OV} which is as follows: the interleaved of a set $X\subset \overline{\mathbb{R}}^n$ is the set of all finite sums $ \sum\limits_{i=1}^{m} e_i\cdot x_i$, with $x_i\in X$ and $\{e_1,\cdots, e_m\}$ a   complete set of mutually orthogonal idempotents in $ \overline{\mathbb{R}}$, i.e., $e_i\cdot e_j=0$ if $i\neq j$ and $\sum \limits_{i=1}^m e_i =1$.   Note that the latter definition is not exactly the one used in Algebra. We stick to this one since in $\ \overline{\mathbb{K}}$ there are no primitive idempotents.  We extend the definition of interleaving  allowing that the number of idempotents involved in the sum of  the  interleaving is countable and not necessarily finite.  Interleavings  can also be done with hypersequences and elements of $\ C^{\infty}(\Omega)$ (see the next section). The expressions  {\it entanglement} or {\it entertwine}  express  the same idea, since several points are connected in the same net that cannot be undone since it is a point in generalized space-time.  This  is actually the way that new points in generalized space-time are created.   Observing a point $\ x$ corresponds to the creation of the point  $x\cdot \alpha^n$, where $\alpha^n$ is defined in the  second paragraph of the next section. Hence, observing is seeing a part of the interleaving $\ x$.  An observation does not change the part of the point that it observes if and only if $\ \alpha^n$ is an idempotent.
 
Consider again $f(x)=x\delta$. The  $halo(0)=B_1(0)$, or the halo of any other point,  contains information of any subset of  $\  \mathbb{R}^n$ via the histories   $\ \mathbb{R}^n\cdot\alpha^r, \ r>0$ or, in general $\mathbb{R}^n\cdot y, \ y\in B_1(0)$ and their  interleavings.  In the same way,  information of any subset of $\mathbb{R}^n$ is contained in the complement of $\ \overline{B}_1(0)$ via histories  $\ \mathbb{R}^n\cdot\alpha^r, \ r<0$.   This can also be seen using the homeomorphisms of $\ \overline{\mathbb{K}}^n$ like those whose existence are proved in  \cite[Theorem 3.3, Theorem 3.4]{JO}. The same holds for subsets of $\ \overline{\mathbb{K}}$.   We can  {\it entertwine  the history} of points $\ x_0\neq x_1\in \mathbb{R}$ using the notion of interleaving: $x= (x_0e_1+x_1e_2)\alpha, e_1+e_2=1$, the latter being idempotents. As a   result,  the measurement, $f(x)$, is also an {\it entertwine}: $f(x)=x_0\rho(x_0)e_1+x_1\rho(x_1)e_2$ which is what is measured in physical reality. We proceed to give an interpretation of this measurement in probabilistic terms.

Consider an entertwine    $\sum\limits_{i}e_i\cdot x_i$.  For each idempotent $\ e\in {\cal{B}}(\overline{\mathbb{R}})$ involved in the sum,    there exists a set $S_e\in {\cal{S}}$, (see \cite[Definition 4.1]{JO}),  such that $e=\chi_{S_e}$  is the characteristic function of $\ S_e$ (see \cite{JO, JJO}).  If  there exists $\ \eta_0 >0$ such that $\ ]0,\eta_0 ]\cap S_e$ is  measurable,  define $\mu(e):=\lim\limits_{\eta\rightarrow 0}\biggl (\frac{1}{\eta}\int\limits_0^{\eta}\chi_{S_e}d\mu\biggr )$.  Since, in an interleaving, the  idempotents involved  form a   complete set of mutually orthogonal idempotents, it follows that $\sum\limits_i\mu(e_i)=1$. Hence    the $\mu(e_i)$'s can be seen as probabilities and, being  countable in number,   there exists $\ i_0$ such that $\mu(e_{i_0})>0$. The interpretation is that whenever  $\mu(e_i)>0$  the measurement at the corresponding point $\ x_i$ is more  likely to be obtained because $f(x_i)$ will appear with the same probability in the resulting measurement (see the example in the previous paragraph).    We say that  the entanglement is a {\it complete entertwine} if $\ \mu(e_i)>0,\ \forall \ i$ and is a {\it perfect entertwine} if $\ \mu(e_i)=\mu(e_j),\ \forall\ i,\ j$. In the latter case the number of idempotents involved must be finite.  Since measurements involve the same generalized function $\ f$,  the whole history of measurement is determent and cannot be changed unless the entanglement is undone.  One can let the $x_i$'s in an interleaving take values in disjoint subsets $X_i \subset \mathbb{R}^n$  letting  the probabilities relate  to the sets $\ X_i$ which, for example,  can be regions in physical reality. In this case, an interleaving can be seen as a function from $\ x: I\longrightarrow \bigcup\limits_i X_i$. For example, rolling  a  dice produces an entertwine  $\ x=\sum\limits_{i=1}^6 e_i\cdot i$, with $X_i=\{ i\}$,  being  perfect only if  the dice is honest.


\section{A Generalized Fixed Point Theorem}

In this section we prove a fixed point theorem which is one more piece of the Generalized Differential Calculus whose  development started  in \cite{OJRO, OJR}.  All the features of this calculus have been extended in \cite{giordano1} to  the context of  Robinson-Colombeau rings of generalized numbers which includes the fields  $\overline{\mathbb{K}}/{\cal{M}}$, where ${\cal{M}}\lhd \overline{\mathbb{K}}$  is maximal (see \cite{JO, V}).  In \cite{giordano1} a  differential calculus is also developed for the latter  rings and the idea of membrane extended.

In the sequel, ideas contained  in   \cite{gar3, giordano} will be used. Let $\widetilde{\mathbb{N}}\subset  \overline{\mathbb{R}}$ be the set of generalized numbers with a representative in $\mathbb{N}^I$ and  $\widetilde{\mathbb{N}\cup\{\infty\}} $  elements of the form $e\cdot +(1-e)\cdot \infty$, $e^2=1$. Another way to view  these elements is to consider   $\ {\cal{E}}_M(\mathbb{N})$ and factor by the ideal $\ {\cal{I}}$ defined in the previous section.  The elements of $\ \widetilde{\mathbb{N}}$ are called {\it hyper natural numbers}.  In the same way one can define the ring of {\it  hyper integers}   $\ \widetilde{\mathbb{Z}}$.  We extend the notation $\alpha^n, n\in \mathbb{N}$, introduced in \cite{JO}, to the case when $n\in \widetilde{\mathbb{N}\cup\{\infty\}} $: if  $n=[(n_{\varepsilon})] $ then $\alpha^n=[(\varepsilon^{n_{\varepsilon}})]$. We can extend this definition  to   $n \in \widetilde{ \mathbb{Z}}$ as long as the set $\ \{ n_{\varepsilon}<0\}$ is bounded, the reason to require  this being obvious.  With this notation,    idempotents are also of the form $\alpha^n$, where, in this case, $n $ consists of a string of $0$'s and $\infty$'s.

A {{\it hypersequence}} is a map $x : \widetilde{\mathbb{N}}\longrightarrow {\cal{G}}(\Omega)$ and is denoted by  $(x_n)$.   If $\ x(\widetilde{N})\subset \overline{\mathbb{K}}$, we say that $(x_n)$  {{\it converges}} to $L\in \overline{\mathbb{K}}$ if given $r>0$, there exists $n_0\in \widetilde{\mathbb{N}}$ such that if  $n>n_0$  then $L-x_n\in V_r(0)$.  Since we are in a Hausdorff space, limits are unique whenever  they exist.  Such a hypersequence $(x_n)$ is a {{\it Cauchy hypersequence}} if given $r>0$ there exists $n_0\in \widetilde{\mathbb{N}}$ such that if  $m, n>n_0$  then $x_m-x_n\in V_r(0)$.  $\overline{\mathbb{K}}$ being a complete metric space, we have:

\begin{lema}
 Let $x=(x_n)$ be a hypersequence. Then $x$ is a convergent sequence if and only if it is a Cauchy sequence if and only if for each $r\in \mathbb{R}_+^*$  there exists $n_0\in \widetilde{\mathbb{N}}$ such that $n>m>n_0$ implies that $x_n-x_m\in V_r(0)$.

\end{lema} 

The sequence $x= (x_n=\frac{1}{n}), n\in \mathbb{N}$ does not converge in $\overline{\mathbb{K}}$ because its elements form a grid of equidistant points. However,   the hypersequence $x= (x_n=\frac{1}{n}), n\in \widetilde{\mathbb{N}}$,  converges to $0\in \overline{\mathbb{K}}$.  In fact, given $r>0$, take $\ \alpha^{-r}\leq n_0=[(\lfloor\varepsilon^{-r}+1.5\rfloor)]\leq \alpha^{-(r+1)}$.  If $n>n_0$,  we have  that $\frac{1}{n}\in V_r(0)$.  We also know  that $\sum\limits_{n\in \mathbb{N}} \frac{1}{n}$ diverges. However if  one sums  over a countabel subset of $\widetilde{\mathbb{N}}$ containing a finite number of elements of  norm one  then the sum converges. Recall that in this setting a series $\sum a_n$ converges if and only if $a_n \longrightarrow 0$.

Let  $r\in\  ]0,1[$ be fixed and suppose that we want $r^n\in V_t(0)$, where $n=[(n_{\varepsilon})]$.  For this to occur one must have $r^{n_{\varepsilon} }<\varepsilon^t$. From this it follows that $n_{\varepsilon}>(\frac{-t}{|\ln(r)|})\cdot \ln(\varepsilon)$. Hence we may take $n_{\varepsilon}=2\cdot \biggl\lfloor\frac{-t}{|\ln(r)|} \cdot \ln(\varepsilon))\biggr\rfloor$. Note that $n<\alpha^{-1}$ and hence is moderate (actually its norm equals $1$). Clearly, for any $m>n$ we have $r^m\in V_r(0)$. This proves that the hypersequence $(r^n)$ converges to $0$.

Any sequence $\hat{x}_0 \ : \mathbb{N}\longrightarrow \mathbb{K}$ defines a map  $\ \hat{x} \ : \widetilde{\mathbb{N}}\longrightarrow \mathbb{K}^{I}$ in the obvious way: $\hat{x}_n=(\varepsilon\longrightarrow \hat{x}_{0n_{\varepsilon}})$.  If $\hat{x} $ is moderate then it defines a hypersequence  $x$.  If there exists $n_0\in \widetilde{\mathbb{N}}$ and $L\in \overline{ \mathbb{K}}_{as}=\mathbb{K}+ \overline{ \mathbb{K}}_0\ $ (see \cite{JO}) such that $n>n_0$ implies $\ x(n)-x(n_0)\in V_1(0)$, then the sequence  $(\hat{x}_{0n})_{n\in \mathbb{N}}$ converges to $L_0\in \mathbb{K}$,  with $L_0\approx L$,  in the classical sense  and the whole history of measuring this convergence is contained in the sentence "$n>n_0$ implies $\ x(n)-x(n_0)\in V_1(0)$". Conversely, if $(x_{0n})_{n\in \mathbb{N}}$  converges to the real number  $L_0$ and for each $\varepsilon$ one choses $n_{\varepsilon}$ minimal such that $n>n_{\varepsilon}$ implies $|\hat{x}_n-L|<\varepsilon$ then $n=[(n_{\varepsilon})]$ must be an element of $\ \widetilde{\mathbb{N}}$ if $x$ were to converge. This shows that a sequence  of measurements can have precision  $\alpha^{k_0}$, for some $k_0$, but not for  all $k\in \mathbb{N}$. It  might be possible to   infer from this if  classical space-time is  discontinuous (unless we declare it continuous and stop measuring beyond $V_1(0))$.    Since it is not at all clear that the convergence of $\ \hat{x}_0$ implies de convergence of $\ x$, one might question the definition of the notion of integral given in \cite{OJRO, OJR, MK}. We shall prove that, at least in  this case, limits do exists and are equal. It  is important that this is true  if we want that  classical theories  also  hold in the Colombeau environment.  If the  classical sequence $(\hat{x}_{0n})$ converges   to $L_0$ and the hypersequence it induces  converges to $L\approx L_0$, but $L\neq L_0$, then this might be the reason for deviations  or strange behavior of systems where $\ L_0$  is used. We shall come back to this again further in this section when we consider the notion of association in ${\cal{G}}(\Omega)$.  As we shall see, the same caution should be taken when dealing with differential equations.

Looking  at classical  continuity of a function $\  f$ at a point $x_0$,  its  history of  measurements produces a function $\delta(\varepsilon)$ such that $x-x_0\in V_{[\delta(\varepsilon)]}(0)$ implies that $f(x)-f(x_0)\in V_1(0)$ with  $V_{[\delta(\varepsilon)]}(0)$  defined accordingly and noting that continuity implies that we can take the class   $\  [\delta(\varepsilon)]\in  \overline{\mathbb{K}}_0$. If $[\delta(\varepsilon)]\in Inv(\overline{\mathbb{K}})$ then there exists $r>0$ such that $V_r(0)\subset  V_{[\delta(\varepsilon)]}(0)$.  Once again, classical theory stops measuring at $V_1(0)$ and one declares a  classical function continuous in $x_0$  once $f(x)-f(x_0)\in V_1(0)$. Again,   scales less then $\alpha$ are not considered because they are not "detected" by the definitions of classical analysis.        Generalized  Differential  Calculus,  links all scales providing a commonage, generalized space-time.

\typeout{Restricting to compact sets, $\delta(\varepsilon)$ has a lower bound $r$, say,  and thus the previous statements becomes $x-x_0\in V_r$ implies $f(x)-f(x_0)\in V_1$. Conversely, the latter statement captures uniform continuity.  and one observes  that we must consider $f$ as acting on nets.  This is exactly what is captured by the notion of point value and  the topology defined in \cite{top1, top2, MK1}.}

Let $f\in {\cal{G}}(\Omega)$ and let  $[K_{\varepsilon}]$ be a membrane (see \cite{OJRO}). To keep things simple we suppose that $K_{\varepsilon}=K$ is compact for all $\varepsilon$ and contained in an open and relatively compact  subset $\Omega_m \subset \Omega$. Note however that the result obtained also holds for a general  membrane.   Given a fixed  $dV\in \{\frac{1}{n},\ \alpha^r \ :\ r\in I,  n\in \widetilde{N}\}$, consider the partition $P$ of norm $dV$ of $\ K$ contained in $\Omega_m$, i.e., for each $\varepsilon\in I $ and $\ dV<\frac{2}{\mu(\Omega_m)}$, $P_{\varepsilon}$ has norm $dV_{\varepsilon}\in \{  \frac{1}{n_{\varepsilon}} ,\ \varepsilon^r\}$, where $\mu(\Omega_m)$ denotes the Lebesgue measure of $\Omega_m$. Since $K$ is compact, there exists $x_0, x_1\in \widetilde{K}$ such that $m=f(x_0)\leq f(x)\leq M=f(x_1), \forall x\in \widetilde{K}$. Let $s(f,dV), S(f,dV)$ and $S(f,dV,*) $ be, respectively the lower and upper Riemann sum and any other starred Riemann sum with this $dV$ as the norm of the partition.  Then $s(f,n)\leq S(f,dV,*)\leq S(f,dV)$ and $|S(f,dV)-s(f,dV)|\leq \alpha_{-N}\cdot \mu(\Omega_m)\cdot dV$, where $N>0$ is such that $\|(\nabla f)_{|_K}\|_{\infty}\leq \alpha^{-N}$. Choosing $dV < \alpha^{N+r+1}$ (this is possible because the hypersequence $(\frac{1}{n})$ converges to zero), we have that  $|S(f,dV)-s(f,dV)|\in V_{r+1}(0)$ and thus $\{ s(f,dV), S(f,dV),S(f,dV,*)\} \subset V_r(\int\limits_Kf(x)dx)$, where $\int\limits_Kf(x)dx)$ is as defined in \cite{MK, OJRO, OJR}.  From this it follows that the classical and generalized  limits are the same and if $r=1$ then we already have that  $s(f,dV), $ and $\ S(f,dV),S(f,dV,*)$ are all associated to  $\int\limits_Kf(x)dx$ and thus, in the classical sense, they are all equal.  In case we take $dV=\frac{1}{n}$ we have a hypersequence and we just proved its convergence. For each $\ s(f,dV),  \ S(f,dV),S(f,dV,*)$ and  $\int\limits_Kf(x)dx)$  there exists an element $c\in \widetilde{K}$, for  each one of  them there exists one such element,  so that it is equal to $\ f(c)\cdot\mu(K)$. If   $\varphi\in {\cal{D}}(\Omega)$  is non-negative, then  $\int\limits_Kf\varphi (x)dx)=f(c)\cdot  \int\limits_K\varphi(x)dx$, for some $c\in \widetilde{K}$. In particular, if $\varphi\in A_0(\Omega)$ (see \cite{AB})  is positive,  then $\int\limits_Kf\varphi (x)dx)=f(c)$.  This is useful when  looking  at the notion of association in ${\cal{G}}(\Omega)$.

Fix $k\in \mathbb{R}_+^*$ and, in  ${\cal{E}}_M(\mathbb{N})$,  define  $\hat{x}(n)(\varepsilon)=\chi_{S(n)}\cdot \varepsilon^k$, where $\chi_{S(n)}$ is the characteristic function of the set $\ S(n)=\{(n_{\varepsilon})^{-1}: \varepsilon \in I\}$. If $n\in \mathbb{N}$ then  $\ S(n)= \{n^{-1}\}$ and $\hat{x}(n)$ would be non-zero only when $\varepsilon=\frac{1}{n}$, having value $\frac{1}{n^k}$.  If these were measurements or observations, and since we cannot make infinitely many measurements or observations, the result would be points of the sequence $(\frac{1}{n^k})_{n\in \mathbb{N}}$ and hence converge to $0$. However,  if we look at the hypersequence $x$ then     $x(n)=e_n\cdot \alpha^k$, where $ e_n=[\chi_{S(n)}]\in{\cal{B}}(\overline{\mathbb{K}})$. Hence $x(n)=0, \forall\ n\in \mathbb{N}$ and $\|x(n)-x(m)\|=e^{-k}\cdot \|e_n-e_m\|\in \{0,e^{-k}\}$. Consequently, this is not a converging hypersequence.  So it might be that measurements in classical space-time are  incomplete.  

In the introduction of \cite{norman} there are two examples that are worth rewriting into this context.  The first is that of a point mass of weight $1$ at a  point $x_0$ on the real axis.  Consider the membrane   $M=V_1(x_0)$. Then the corresponding functionals can be written as 

$$L(\varphi)=\frac{1}{vol(M)}\int\limits_M\varphi(x)dx=2\alpha^{-1}\int\limits_M\varphi(x)dx$$

\noindent By \cite[Proposition 3]{OJRO} we have that $L(\varphi)=\varphi(c)$, for  some $c\in M$. Since $\varphi\in {\cal{D}}(\mathbb{R})$, we have that $\ \varphi (M)\subset B_1(\varphi(x_0))=halo(\varphi(x_0))\ $ and hence $\ \varphi (c)\approx \varphi (x_0)$, the latter being   what measurements  give us,  but   $\ \varphi(c)\in\overline{\mathbb{R}}$ being the actual value.  We can extend this to $\ \varphi \in {\cal{G}}(\mathbb{R})_{as}$ (which will be defined yet in this section)  by  taking  $\ M=V_k(x_0)$ such that $\ \varphi (M)\subset B_1(\varphi(x_0))$.

The other example is that  of  a dipole at $0$  with moment $1$. In this case,  the functionals can be written as 

$$ M(\varphi)= \alpha^{-1}\cdot(\varphi(\alpha)-\varphi(0))$$
 
\noindent This amounts to $\varphi (\alpha)=\varphi(0)+M(\varphi)\cdot \alpha$. Expanding around $0$ gives us that $M(\varphi)-\varphi^{\prime}(0)\in V_1(0)$, i.e., $M(\varphi)\approx \varphi^{\prime}(0)$,  the latter being   what we measure,   the actual value being  $\ M(\varphi)\in \overline{\mathbb{R}}$. Again, this   can  be  extended   to ${\cal{G}}(\mathbb{R})$. In both cases,  our measurements are the only "real" point  in the halo of the actual values.   Howbeit,   classically we cannot differentiate among points in a halo. From the classical point of view, each halo contains only one point but, as we already saw, they do interfere in physical reality.  In both examples,  the measurable effect in physical reality  is a product of an infinity, $2\alpha^{-1}$ respectively  $\alpha^{-1}$, and an infinitesimal, $\int\limits_M\varphi(x)dx$ respectively  $\varphi(\alpha)-\varphi(0)$, both accredited   and coexisting in harmony in this generalized milieu.  

In the generalized  environments   the history of  measurements and of convergence is captured  and not each measurement  separately (however arbitrarily).   If  the problem lies in the classical  mathematical tools  that we use and not in  space-time its self,    then we hope that these examples  help to convince   that  the Colombeau environments, in particular generalized space-time,   are  environments that maybe should be considered. See also \cite[Section 1.6]{MK} and the references mentioned there.

For the reader's sake, we recall the basics  of the sharp topology (see  \cite{scar, hebe, top1, top2, AGJ}).  Let $(\Omega_m)$ be an exhaustion of relatively compact and open subsets  of $\Omega\subset \mathbb{R}^n$. Given $\ f\in {\cal{G}}(\Omega)$ and $(m,p)\in\mathbb{N}^2$,  define $V_{mp}:=Sup\{a\in \mathbb{R}\ |\ \forall \beta\in  \mathbb{N}^n, \ |\beta|\leq p, \|\partial^{\beta}f(\varepsilon, \cdot\ )\|_{\Omega_m}=o(\varepsilon^a), \ \mbox{for }\ \varepsilon\  \mbox{small}\ \}$ and $D_{mp}(f,g):=exp(-V_{mp}(\hat{f}-\hat{g}))$, where $\hat{ }$ stand for representative. The latter are pseudo-ultrametrics  defining the Biagioni-Scarpal\'ezos  sharp topology on  $\ {\cal{G}}(\Omega)$ (see   \cite[Definition 1.9, Proposition 1.10 and 1.11]{JA1}).  This topology is  proved to be equivalent to the topologies given in \cite{top1, top2, AGJ}.  As observed in these references, this topology is metrizible and, with the notation  given above, an ultrametric in $\  {\cal{G}}(\Omega)$  is given by

 $$D(f,g)=sup\biggl\{\ \frac{2\cdot D_{mm}(f,g)}{1+D_{mm}(f,g)} \ :\ m\in\mathbb{N}\biggr\}.$$ 
 
  With  the notation of  \cite{top1, top2, AGJ}, $W^k_{m,r}(0)=\{f\in {\cal{G}}(\Omega)\ |\ | \partial^{\beta}f(x)|\in V_r(0),\ \forall\  |\beta|\leq k,\ \forall \ x\in \widetilde{\Omega_m}\}$,  $V_r(0)=\{x\in\overline{\mathbb{K}}\ :\ |x| <\alpha^r\}$  and $\ \|\partial^{\beta}f\|_{\beta,m}:=[\varepsilon\longrightarrow \| (\partial^{\beta}f_{\varepsilon})_{ |_{\Omega_m}}\|_{\infty}]$.   By \cite[Theorem 3.6]{top2}, the  sets $W^{k}_{m, r}(0)$  form a filtered basis of the sharp topology of $\ {\cal{G}}(\Omega)$.

    If $f, g\in {\cal{G}}(\Omega)$ then $f\approx g$ iff $\ \int\limits_{\Omega}(f-g)\varphi dx\in \overline{\mathbb{K}}_0$ and $f\sim g$
   iff $\ \ \int\limits_{\Omega}(f-g)\varphi dx=0\in \overline{\mathbb{K}}$, $ \forall \varphi \in{\cal{D}}(\Omega)$ (see \cite{MK}). In the latter case,  one says that $f$ and $g $ are {\it test equivalent} and in the former case one  says that they are {\it associated}.   In \cite{JO},     $\ \overline{\mathbb{K}}_0=\{x\in\overline{\mathbb{K}}\ :\ x\approx 0\}$ was first formally  given a notation as was $\overline{\mathbb{K}}_{as}=\mathbb{K}+\overline{\mathbb{K}}_0$  (see \cite[Definition 1.6.5]{MK} where they were first defined). Here, we  introduce the notation   $\  {\cal{G}}(\Omega)_0=\{f\in  {\cal{G}}(\Omega) \ :\ f\approx 0\}$  and $\  {\cal{G}} (\Omega)_{as}=\{f\in  {\cal{G}}(\Omega) \ : \ \exists\  T\in {\cal{D}}^{\prime}(\Omega),\ \mbox{such that}\  f\approx T\}$. Clearly,    $\  \overline{\mathbb{K}}_0\subset {\cal{G}} (\Omega)_0$.  Define the {\it halo} of $\ f\in {\cal{G}}(\Omega)$ as $halo(f):=f+B_1(0)=B_1(f)$, where $B_1(0)$ is   the ball $\ \{f\in {\cal{G}}(\Omega)\ :\ D(f,0)<1\} \subset  \ {\cal{G}}(\Omega)$ (see also  \cite{jailson}). Confusion should not arise between  $\  {\cal{G}}(\Omega)_0$  and    ${\cal{G}}_0(\Omega)$, the original Colombeau algebra (see \cite{MK}).  Association in Colombeau algebras has  been presented as an algebraic notion substituting equality in some sense.  We proceed to  prove that it is in fact  a topological notion  and    use this to prove  that, in  a topological sense, Schwartz generalized functions, and hence classical  solutions of differential equations,  are scarce.  
  
  
  \begin{prop}

Let  $f, g\in {\cal{G}}(\Omega)$.  Then the following hold.

\begin{enumerate}
\item   If $\ g\in halo(f)$ then $g\approx f$.
\item $B_1(0)=halo(0)\subset {\cal{G}}(\Omega)_0$.
\item If $\ \| f\|_m \in \overline{\mathbb{K}}_0,\ \forall\ m$, then $f\in {\cal{G}}(\Omega)_0$.
\item If   $\ Im(f) \subset  \overline{\mathbb{K}}_0$ then $f\in  {\cal{G}}(\Omega)_0$.
\item ${\cal{G}}(\Omega)_{as} = {\cal{D}}^{\prime}(\Omega)+{\cal{G}}(\Omega)_0$. 

 \end{enumerate}

\end{prop}  

\begin{proof}    To  prove the first two items,    by hypothesis,  $h:=f-g\in B_1(0)$.   Since $D(h,0)<1$, it  follows that $D_{mm}(h,0)<1$, $\forall m\in \mathbb{N}$ and hence there exists a decreasing sequence $((r_m)\ ,  r_m >0, \ \  \forall\  m )$,   such that $\ V_{mm}(h)=r_m>0, \forall m\in \mathbb{N}$.  Consequently,  $\ h\in \bigcap\limits_{m\in \mathbb{N}, }W^m_{m,r_m}$.  In particular, $h(x)\in V_{r_m}(0), \forall x\in \widetilde{\Omega_m}$, i.e., $|h(x)|< \alpha^{r_m}, \forall x\in \widetilde{\Omega_m}$. Let $\varphi\in {\cal{D}}(\Omega)$ and $m\in \mathbb{N}$  be such that $ supp(\varphi)\subset \Omega_m$;  then  $|\int\limits_{\Omega}h\varphi dx| \ =|\int\limits_{\Omega_m}h\varphi dx| \leq \alpha^{r_m}\cdot \|\varphi\|_{\infty}\cdot \mu(\Omega_m)$.  Hence $\int\limits_{\Omega}h\varphi dx\in \overline{\mathbb{K}}_0$.

If $\varphi\in {\cal{D}}(\Omega)$ then there exists $m$ such that $supp(\varphi)\subset \Omega_m$. Hence \\ $|\int\limits_{\Omega}f(x)\varphi(x)dx|= |\int\limits_{\Omega_m}f(x)\varphi(x)dx|\leq \|f\|_m\cdot \|\varphi\|_{\infty}\cdot\mu(\Omega)\in \overline{\mathbb{K}}_0$, proving the third item.

To prove the fourth  item, use an appropriate Riemann sum, as was shown to exist  in the examples preceding the proposition, and the fact that ${\overline{\mathbb{K}}_0}$ is a ring. It also follows by the previous item. The fifth item is an obvious statement.
 \end{proof}

\typeout{{\bf daqui para baixo ainda nao est\'a certo.}
  To prove the third item, suppose that $f\not =0$ in $\Omega$, then we may suppose that $supp(f)$ is compact and that  there exists $x_0\in \Omega$ such that $f(x_0)\neq 0 \in \overline{\mathbb{R}}$.  Given a mollifier   $\rho$  (or any element in $ A_0(\Omega)$), we have that    $\ \int\limits_{\Omega}f(x)\rho_{\varepsilon}(x-x_0)dx=0, \ \forall \ \    \varepsilon\in I$.   On the other hand, choosing $m>1 $ such that $supp(f)\subset \Omega_m$ and $x_0\in \widetilde{\Omega_m}$, we have that $\ \int\limits_{\Omega}f(x)\rho_{\varepsilon}(x-x_0)dx-f(x_0)= \ \int\limits_{\Omega}(f(x_0+\varepsilon z)-f(x_0))\rho(z)dz\leq   \int\limits_{\Omega}|\langle \nabla f(z_{\varepsilon})|\varepsilon z\rangle|\cdot |\rho(z)|dz\leq         \sqrt{n}\cdot   \alpha^{\frac{r_m}{2}}\cdot diameter(\Omega_m)\cdot\|\rho\|_{\mathbb{L}_1} $. Since $r_m>0$ we have that $f(x_0)=0$, a contradiction.  
    a
   To prove the fourth item, take $\varepsilon =\frac{1}{k}$ in the preceding paragraph obtaining a sequence $z_k:=   \int\limits_{\Omega}f(x)\rho_{\frac{1}{k}}(x)dx \in \overline{\mathbb{K}}_0$ converging to $f(x_0)$. Note that $z\in \overline{\mathbb{K}}_0$ if and only if $z<\frac{1}{k},\ \forall \ k\in \mathbb{N}$.  Since $z_n \longrightarrow z:=f(x_0)$, for each $q>0$ there exists $n_q\in \mathbb{N}$ such that $n>n_q$  implies that $z-z_n\in V_q(0)$. So given any  $l\in \mathbb{N}$,  it follows that $\ |z|<|z_n|+\alpha^q< \frac{1}{l}+\alpha^q$ and hence $|z|\leq \frac{1}{l}$. Consequently, $z\in \overline{\mathbb{K}}_0$.
   a
   The last item follows by the fourth because $\overline{\mathbb{K}}_0$ is a subring of $\ \overline{\mathbb{K}}$.
a
If $f(x)=x\delta\in {\cal{G}}(\mathbb{R})$ then $f(x)=0\forall\  x\in \mathbb{R}$ as can be seen easily (see \cite{MK}). However $f(\alpha)=\rho(1)$, where $\rho$ is the mollifier defining $\delta$.
a
{\bf verificar o que est\'a acima.}}

  If $f\approx 0$ then for each $x_0\in \Omega$  and each $B_r(x_0)\subset \Omega$ we have that there exists $c_r\in \widetilde{B_r(x_0)}$ such that $f(c_r)\in \overline{\mathbb{K}}_0$.  This follows   from the observation  at the end of the paragraph about  Riemann sums.
  
  In \cite{OJR} it was proved that $\ \mathbb{R}^n$ is a discrete subset of $\ \overline{\mathbb{R}}^n$ and that if $r\in\mathbb{K}^*$ and $x\in \overline{\mathbb{K}}$ then $\| rx\|=\| x\|$.  Our next results state  that the same  is true for $D^{\prime}(\Omega)$ and $C^{\infty}(\Omega)$.    The first part  of the next corollary  is  in fact nothing more than a topological interpretation of \cite[Proposition 1.6.3]{MK} and \cite[Proposition 1.7.28]{MK}. The second part looks at the building blocks of $\ {\cal{G}}(\Omega)$ and equate  them with the building blocks of $\ \overline{\mathbb{K}}$.

  \begin{cor}

Let $\Omega\subset \mathbb{R}^n$ be an open subset. The embedding of  $\ D^{\prime}(\Omega)$ in  $\ {\cal{G}}(\Omega)$ is a discrete embedding. Moreover, if   $\ h\in C^{\infty}(\Omega)^*$ and $f\in {\cal{G}}(\Omega) $  then $\|hf\| = \|f\|$.

\end{cor}

\begin{proof}  Since the embedding is linear, we just have to prove discreteness at the origin.   Let $f\in B_1(0)\cap{\cal{ D}}^{\prime}(\Omega)$.  By the previous proposition,   $f\approx 0$ and hence,  by \cite[Proposition 1.6.3]{MK},  $f=0$.   The proof of the second part is straightforward. \end{proof}

Since   $g\approx f$ if   $g\in halo(f)$,  it follows that  uniqueness of solutions and association  can be seen as  statements about halos. That is, the weaker notion of equality called association  is a topological notion!  Equality is an algebraic notion! Elements of ${\cal{G}}(\Omega)$  that are associated are  indistinguishable  one from another from the  classical point of view  (See also \cite[Section 1.6]{MK}).         An element  of  $\ {\cal{G}}(\Omega)$, not in the halo  of any point of  $\  {\cal{D}}^{\prime}(\Omega)$,  is  called a vampier in \cite{MK}: it has no distributional shadow. Note however that it can  be that, multiplying it with an $\ e\in {\cal{B}}(\overline{\mathbb{R}})$ it has a shadow, even infinitely many and thus, it can flaunt  omnipresence (see the notion of  entertwine  or of  support  in the next section).    This shows that the construction of $\ \widetilde{\Omega}_c$ starting from $\Omega$ is the same as the construction of  $\  {\cal{G}}(\Omega)$ starting from $\ {\cal{D}}(\Omega)$. All  that is classical becomes discrete in the  generalized environments.  Once again, no sequence converging  in $D^{\prime}(\Omega)$ can converge in ${\cal{G}}(\Omega)$ and, classically,  measuring convergence is stopped at the sets $\ W^m_{p,1}$!  Again, it appears that  one has to call upon  hypersequences  to fix this.

That being so,  even in the Colombeau Algebras there is just one notion of equality, i.e.,  the classical one!  Association is not equality in any sense but a topological statement and test equivalent looks more like classical  equality but is not. In fact, once again,  let $f=x\delta$, $ q\in\mathbb{N}$  and let the delta-net be induced by the mollifier $\ \rho$.  Then $\int\limits_{\Omega}f_{\varepsilon}\varphi dx=\varepsilon\cdot \int\limits_{\Omega}z\rho(z)\varphi(\varepsilon z)dz$ $= \varepsilon\cdot \int\limits_{\Omega}z\rho(z)[\varphi(\varepsilon z)-\varphi(0)]dz$ $= \frac{\varepsilon^q}{(q-1)!}\cdot \int\limits_{\Omega}z\rho(z)\varphi(\theta\cdot \varepsilon z)dz=o(\varepsilon^q)$.  This proves that $\int_{\Omega}f\varphi dx=0\in \overline{\mathbb{K}}$ and thus, $f\sim 0$ but $\ f\neq 0$.  For this $f$, $f(x)=0,\ \forall \ x\in \mathbb{R}$ but  it is not true that it is zero uniformly on compact subset of $\ \mathbb{R}$, i.e., it is not negligible  at level zero and even more, $\ f\notin W^0_{m,r}(0)$ with $\  r>0$.    Showing that, classical mathematical tools are incomplete  not being able to measure    the full physical reality of space-time.  
 
There should be no strangeness in the fact that  $\ f(x)=x\delta$ is identically zero on $\ \mathbb{R}$ but is not in $\widetilde{\mathbb{R}}_c$.  Recall  that $ \ \mathbb{R}$ forms  a grid of equidistant points in $\ \overline{\mathbb{R}}$ and thus $\ f=0$ on a discrete grid with no   accumulation points.  However,  $f^{\prime}(x)=\delta+x\delta^{\prime}$ is not necessarily   identically zero in $\ \mathbb{R}$ because $\ f^{\prime}(0)=\rho(0)\cdot\alpha^{-1}$.  For comparison,  $g(x)=sin(\pi x),\ x\in \mathbb{Z}$, is identically zero but  $\  g^{\prime}=\pi cos(\pi x),  \ x\in \mathbb{Z}$,   is not.  The most  puissant   differential calculus   developed in $\ \mathbb{Z}$ will not be able to detect that    $\ g$  is not identically zero but  that its derivatives effect  reality in $\ \mathbb{Z}$.

 In case of $\ \overline{\mathbb{K}}$,  nets of $\mathbb{K}$  are its  building blocks   and one knows  that $\ \mathbb{K}$ embeds  as a grid of equidistant point into  $\ \overline{\mathbb{K}}$. In case of $\ {\cal{G}}(\Omega)$, the building blocks are nets of elements of $\ C^{\infty}(\Omega)$ (see also \cite{AB, hebe, JFC2, MK}  about how to construct intrinsically $D^{\prime}(\Omega)$ starting with the embedding of $C^{\infty}(\Omega)$). Hence  the following   results should not come as a surprise (see also \cite[Theorem 5.8]{JJO}).

 \begin{cor}
 The elements of $\ C^{\infty}(\Omega)$ form a grid of equidistant points in $\ {\cal{G}}(\Omega)$. In particular, ${\cal{G}}(\Omega)$ is a fractal.
 \end{cor}
 
 \begin{proof} Since $\ C^{\infty}(\Omega)$ is diagonally embedded it follows that $V_{mp} = 0, \ \forall m,p$.  By  \cite[Lemma 3.6]{JO}, it follows that the inductive dimension  $dim({\cal{G}}(\Omega))=\infty $ and, being an ultrametric  space,     $Ind({\cal{G}}(\Omega))=0$. Consequently, $\ {\cal{G}}(\Omega)$ is a fractal.   \end{proof}

 Generalized Differential Calculus  aims at  relying  soly  on intrinsic results and not on classical results to prove existence and uniqueness of differential equations, among other things. Should it be that  results proved   using classical tools to prove moderateness and establishing solutions for differential equations should be considered incomplete?  Our next endeavor  is to show how to sidestep this, at least in the direction that   we are going to investigate,  in such a way that it  extends   naturally   all classical results in this direction. Might it  be the way to complete  results obtained until now using classical tools?

Given a net of maps $(T_{\varepsilon})$ and $n=[(n_{\varepsilon})]\in \widetilde{\mathbb{N}}$, we denote by $T^n$ the net $(T_{\varepsilon}^{n_{\varepsilon}})$ acting on nets $f=(f_{\varepsilon})$ as  $\ T^n(f)=(T^{n_{\varepsilon}}_{\varepsilon}(f_{\varepsilon}))$.    Suppose that $T^n$   is  well defined in ${\cal{G}}(\Omega)$, with $\Omega\subset \mathbb{R}^n$ an open subset,  and denote  by $T$ the net when $n=1$. Let ${\cal{A}}\subset {\cal{G}}(\Omega)$ be a closed subspace and suppose that the restriction $\ T_{|_{\cal{A}}} : \cal{A}\longrightarrow \cal{A}$.  The map   $\ T :  \cal{A}\longrightarrow \cal{A} $ is said to be a {{\it  contraction}} if there exist $L\in \overline{\mathbb{R}}_+^*, \lambda   \in ]0,1[\ $ such that $\ L< \lambda  $   and  $\ |T(f)-T(g)|(x)\leq L\cdot|f-g|(x)$.  It easily follows that  $\ |T^n (f)-T^n(g)|(x)\leq L^n\cdot|f-g|(x)\leq {\lambda}^n\cdot|f-g|(x)$. Our  interest is  when the      hypersequence  $(T^n(f))$ converges in $\ {\cal{G}}(\Omega)$.    We  look at some  examples  inspired by the classical analog.

 Let  $\Omega\subset \mathbb{R}^n$ be an open subset, $L\in \mathbb{R}_+$  and $x_0\in  \widetilde{\Omega}_c$. Define ${\cal{A}}=\{ f\in {\cal{G}}(\Omega)\ : \ |f(x)-x_0|\leq L, \ x\in \widetilde{\Omega}_c\}$ and consider it with the induced sharp topology. Then ${\cal{A}}$ is a closed subset of  $ {\cal{G}}(\Omega)$. In fact, let  $(f_n)$ be a Cauchy sequence in ${\cal{A}}$.  Since ${\cal{G}}(\Omega)$ is complete (see \cite{top1, top2, AGJ}),    $f_n \longrightarrow f$, for  some $\ f$,  and thus    $x_n:=|f_n(x)-x_0|\longrightarrow |f(x)-x_0|=:a$.  Since $x_n$ converges to $\ a$, for each $r>0$ there exists $n_0$ such that if $n>n_0$ then $a-x_n\in V_r(0)$ and thus $|a-x_n| < \alpha^r$. Hence $a<x_n +\alpha^r \leq L+\alpha^r$. It follows that $a\leq L$, i.e.,  $\ f\in {\cal{A}}$, thus proving that ${\cal{A}}$ is a closed subset  of $\ {\cal{G}}(\Omega)$.  We did not use the fact  that $L\in \mathbb{R}$ but the reason why it appears here is that when considering some differential equations,  compositions of generalized  maps is necessary. In \cite{OJR} it was shown that this results in the classical composition of maps. Since domains of generalized functions are $\ \widetilde{\Omega}_c$, the real bound $\ L$ will guarantee   that   compositions of maps  are allowed.   Consequently, if $f\in {\cal{A}}$ then $|f(x)|\leq L+|x_0|$ proving that $\| f\| \leq 1$, i.e.,   ${\cal{A}}\subset \overline{B}_1(0)$. 
 
 Suppose that for each $\varepsilon \in I$ we have a map  $T_{\varepsilon}\ :\  {\cal{A}}_{\varepsilon} \longrightarrow {\cal{A}}_{\varepsilon}$, with ${\cal{A}}_{\varepsilon}= \{h\in C(\Omega,\mathbb{R})\ :\  |h(x)-x_{0\varepsilon}|<L\} $,  which   is a contraction with Lipschitz constant $K_{\varepsilon}<\lambda \in\  ]0,1[$, the latter being independent of $\varepsilon$. It is clear that with these settings $T^n$, with $n\in \widetilde{N}$,  is a  well defined  and continuous map  in $ {\cal{G}}(\Omega)$ and  $T^{n+m}=T^n\circ T^m, n,m\in \widetilde{N}$. In particular, $T^{n+1}=T\circ T^n$, observing that $1$ must be seen as an element of $\ \widetilde{N}$.

Another classical situation is the following.  The set  $\ \Omega\subset \mathbb{R}^n\ $ is open and relatively compact, for  each $\varepsilon\in I$,  $T_{\varepsilon}(x)(t)=x_{0\varepsilon}+\int\limits_{t_0\varepsilon}^t h(s,x(s))ds, t\in \Omega$, $h\in C^{\infty}(\mathbb{R}^{n+1})$  and ${\cal{A}}_{\varepsilon}=C(\Omega, \mathbb{R}^n)$ (see \cite{Erl1, Erl2} and \cite[Theorem 3.1]{erlacher}). In this case, there exists $n_{0\varepsilon} \in \mathbb{N}$,  such that $T_{\varepsilon}^{n_{0\varepsilon}}$ is a contraction  with Lipschitz constant $K_{\varepsilon}\leq \lambda <1$, the latter being real and  fixed. If $n_0:=[(n_{0\varepsilon})]$, then  define $T=[(T_{\varepsilon}^{n_{0\varepsilon}})]$ and thus reducing it  to the case of   the previous paragraph,   a classical argument.

\begin{teo}[{\bf{The Generalized Fixed point Theorem}}]

Let  $\Omega \subset \mathbb{R}^N$,  ${\cal{A}}\subset {\cal{G}}(\Omega)\cap \overline{B}_1(0)$ a close subset.  For each $\varepsilon \in I$, let  ${\cal{A}}_{\varepsilon}= C^{\infty}(\Omega)\cap {\cal{A}}$, initial conditions taken for that specific $\varepsilon$,  and  $\ (T_{\varepsilon} )$  a net of functionals from  ${\cal{A}}_{\varepsilon}$ to its self.  If there exists $k\in \widetilde{\mathbb{N}}$ such that each $\ T_{\varepsilon} ^{k_{\varepsilon}}$ is Lipschitz with Lipschitz constant $K_{\varepsilon}<\lambda\in \ ]0,1[$,  then $T=\ (T_{\varepsilon}^{k_{\varepsilon}} )$  is well defined, continuous  and   has a unique fixed point $\tilde{x}\in {\cal{A}}$. 
\end{teo}

\begin{proof} The proof  uses  what was  already  discussed and also  freely facts    about the topology.
 
Choose   any $x\in {\cal{A}}$ and consider the hypersequence $(T^n(x))$. Given $ r>0$ and $m\in \mathbb{N}$,  choose $n_0\in \widetilde{\mathbb{N}}$ such  that $\ {\lambda}^{n_0}\cdot \alpha^{-1} \in V_{4^{mN}r}(0)=V_{4^Nr_1}(0)$, where  $r_1=4^{N(m-1)}r$. If  $n, s>n_0$ then, writing  $n=n_0 +k,\  s=n_0+l$, we have that $|T^n(x)-T^s(x)|(t) \leq$ $ {\lambda}^{n_0}\cdot |T^k(x)-T^l(x)|(t)\leq  {\lambda}^{n_0}\cdot \alpha^{-1}\in V_{4^Nr_1}(0)$.  Setting  $F(t)= (T^k(x)-T^l(x))(t)$, with $\ t\in \Omega_m$,  we have  that  $F\in \overline{B}_1(0)\cap W^0_{m,4^Nr_1}(0)$ and hence $|F^{\prime\prime}(t) |\leq \alpha^{-0.5r_1}$.   Without loss of generality, we may   considering $N=1$, obtaining, since $\widetilde{\Omega_m}$ is open,    $F(t+\alpha^{2r_1})-F(t)=F^{\prime}(t)\cdot\alpha^{2r_1}+\frac{ F^{\prime\prime}(c)}{2}\cdot\alpha^{4r_1}$ and thus  $\ |F^{\prime}(t)|=$  $ \biggl |\frac{F(t+\alpha^{2r_1})-F(t)}{\alpha^{2r_1}}+\frac{ F^{\prime\prime}(c)}{-2}\cdot\alpha^{2r_1}\biggr |\leq \frac{2\alpha^{4r_1}}{\alpha^{2r_1}}+\frac{\alpha^{-0.5r_1}}{2}\cdot \alpha^{2r_1}<\alpha^{r_1}$. This proves that $F\in W^1_{m,r_1}(0)=W^1_{m, 4^{m-1}r}(0)$. By induction, we have  that   $F\in W^m_{m,r}(0)$.  This proves that $\ (T^n(x))$ is a Cauchy hypersequence and hence converges. Since $ \ T^{n+1}=T\circ T^n$ and $T$ is continuous,  the limit is a fixed point of $\ T$.
\end{proof}

In case   compositions of maps are  not involved,  the theorem still holds if ${\cal{A}}\subset B_R(0),\ R>1$ (It  always holds since estimates are made in $\ \widetilde{\Omega_m}$).   Together with the other tools, it  makes  the Generalized Differential Calculus developed in \cite{OJRO, OJR} (see also \cite{spjornal})   a useful tool  to generalize most classical results.   For example,  one can prove the local, and hence the global, existence  of geodesics in $\widetilde{M}_c$ (see the next section).  If the solution, given by the theorem,  is not a classical solution but it has  a distributional shadow, it is a non-vampire,  this might lead one to change $=$ in the equation  by $\approx$. Actually, this is often done but maybe it should be rethought.   Why not  solve the equation as it is in the generalized  environment and then see if it has a distributional shadow?  The distributional shadow will  not be the solution but will behave like the non-vampire.  Having observed   this,    there exist  classical equations which do not have a classical solution but do have a generalized solution which is a non-vampire or omnipresent.    There are  eidolons roaming  under and over  classical radars,  verbalized by  the ancients and Dirac   as atoms and infinities, which,  at rendezvous,  tellingly   influence physical reality.  

 Let's  formalize the argument used in the last part of the  proof of the theorem since it is a useful tool to  be considered in Generalized Analysis.    

\begin{defini} 
 {\bf{The Down Sequencing Argument}} \\ 
 Let $f\in {\cal{G}}(\Omega)$, with $\Omega\subset \mathbb{R}^n$.  If    $\ f\in W^0_{m,r}(0)$ with $r>0$ and  $p_0\in \mathbb{N}^n$  then  $\ f\in W^{|p_0|}_{m,s}(0)$, where $s=4^{-n|p_0|}\cdot r$, i.e.,     $W^0_{m,r}(0)\subset W^{|p_0|}_{m,s}(0)$.
\end{defini}

Using the {\bf{DSA}},     another proof  of a fact already mentioned can be given: If $\ f$ is moderate and negligible at level $0$, then $\ f$ is negligible. This can be considered a statement about rigidity.  In fact,  such an $f\in W^0_{m,r},\ \ \forall\  m, \ r>0$ and hence DSA gives   what was  claimed. In other words   $\ \bigcap\limits_{r>0} W^0_{m,r}=\{0\}$,  showing that these  sets serve as a basis for the sharp topology.     Let's consider the   example that inspired the Generalized Fixed Point  Theorem. Consider the following   equation from \cite{MK2, MK3}.

$$\ddot{x}(t)=f(x(t))\delta(t)+h(t)$$
$$x(-1)=x_0$$ 
$$\dot{x}(-1)=\dot{x}_0$$

\noindent with $h, f\in C^{\infty}(\mathbb{R})$ and $\ \delta$ the Dirac function.  This is a typical differential  equations  from Physics  and Engineering  having a product of distributions in its data and does  not allow the use of classical tools to obtain a solution.      In the references just mentioned, one   has to prove moderateness and other nontrivial steps  involving non trivial classical results. Would this make  such proofs incomplete?   We show that proving moderateness is  not needed and thus suggesting a possible remedy.    Let $b>0$,  $C>0$ a positive  constant   limiting the $\mathbb{L}_1$ norm of the $\ delta-$net, $M:=\int\limits_{-2}^1\int\limits_{-2}^1 |h(r)|drds$,  $L:= b+M+|\dot{x}_0| +\| f\|\cdot C $  and $1+a=min\{\frac{b}{C\cdot\|f\|+|\dot{x}_0|}, \frac{1}{2CK}, 2\}$, where $K$ is a Lipschitz constant of $\ f$ on a compact subset of $\ \mathbb{R}$ containing $\Omega =\ ]-1-\frac{a}{2},\frac{a}{2}[  $. The norm $\| f\|$ is also taken over the same  compact subset.  Let ${\cal{A}}=\{x=(x_1,\cdots,x_n)\in C^{\infty}(\widetilde{\Omega}_c, \tilde{\mathbb{R}}^n)   \ :\ x_i\in {\cal{G}}(\Omega)\ \mbox{and} \ |x(t)-x_0|\leq L\}$. Since   $C^{\infty}(\widetilde{\Omega}_c, \tilde{\mathbb{R}}^n)\cong (C^{\infty}(\widetilde{\Omega}_c, \tilde{\mathbb{R}}))^n$ we have that ${\cal{A}}\subset ( {\cal{G}}(\Omega))^n$ and the latter is a complete metric space (see \cite{top1, top2, AGJ}). From what we already discussed, it follows that ${\cal{A}}$ is a closed subspace of a complete algebra. Define ${\cal{A}}_{\varepsilon}$ accordingly  and let $T_{\varepsilon}$ be defined in ${\cal{A}}_{\varepsilon}$ by  $T_{\varepsilon}(x)(t)=x_{0\varepsilon}+\dot{x}_0(t+1)+\int_{-1}^t\int_{-1}^s f(x(r))\rho_{\varepsilon}(r)drds+\int\limits_{-1}^1\int\limits_{-1}^s h(r)drds$.  We have that   $|T_{\varepsilon}(x)(t)-x_{0\varepsilon}| =| \dot{x}_{0\varepsilon}(t+1)+\int_{-1}^t\int_{-1}^s f(x(r))\rho_{\varepsilon}(r)drds+\int\limits_{-1}^t\int\limits_{-1}^s h(r)drds|$ $\leq | \dot{x}_{0\varepsilon}|\cdot(|t|+1)$  $+\int_{-1}^t\int_{-1}^s| f(x(r))\rho_{\varepsilon}(r)|drds$ $+\int\limits_{-1}^t\int\limits_{-1}^s |h(r)|drds$ $\leq |\dot{x}_{0\varepsilon}|\cdot(2+a)+(2+a)\cdot \| f\|\cdot C+M$ $= |\dot{x}_{0\varepsilon}|\ +M +(1+a)(|\dot{x}_{0\varepsilon}|+\| f\|\cdot C)$ $\leq |\dot{x}_{0\varepsilon}|\ +M +b+\| f\|\cdot C) =L$. It is also Lipschitz:   $|T_{\varepsilon}(x)(t)- T_{\varepsilon}(\tilde{x})(t)|=| \int_{-1}^t\int_{-1}^s (f(x(r))-f(\tilde{x}(t)))\rho_{\varepsilon}(t)drds  |$ $\leq \| f\|\cdot KC\cdot(2+a)\cdot |x(t)-\tilde{x}(t)|  $ $\leq 2 \| f\|C\cdot(1+a)\cdot K |x(t)-\tilde{x}(t)|\leq  K |x(t)-\tilde{x}(t)|$. It follows that the $T$ they define is    Lipschitz  and hence  well defined and continuous  in ${\cal{G}}(\Omega)$. The first part shows  that $T$ maps ${\cal{A}}$ into itself.  Hence we are in the setting of the Generalized  Fixed Point Theorem.  We have a Lipschitz map and the corresponding hypersequence has a fixed point  which is a solution for the system.  Once the Generalized Fixed Point Theorem is applied, there is no need to prove moderateness. The result follows forthwith!     From the equation it  readily allows that $\  \ddot{x}\approx h+f(x(0))\cdot \delta$, i.e., it is a non-vampier. 
 
Supposing $h=0$,  expand $\ \delta$  at collision time $t=0$  and chose  the mollifier such that $\rho(0)=1$. The equation to be solved is
 \\

$$\ddot{x}(t)=f(x(t))\alpha^{-1}$$
$$x(-1)=x_0$$
$$\dot{x}(-1)=\dot{x}_0$$ 

Its solution is obtained  just  as in  the classical case but using Generalized Differential Calculus. Albeit,  the solution might not be a Colombeau generalized function since it might only be defined  on a membrane. This happens because    the Theory of Generalized Differential Calculus strictly  contains the Theory of  Colombeau Generalized Functions (see \cite{OJRO} for an explicit example).  There exists a  non-constant functions with zero derivate (see \cite{OJR}) behaving like a quanta.  For applications, a solution in this more generalized environment   does provide a solution since  $\varepsilon$ small is enough;  one does not  need an  explicit solution in the Colombeau algebra.   The environment of    the full algebra is  treated in \cite{paulo} relying  on results from (\cite{top1, top2, JRJ, CGS, Ni}). In the case of the invariant algebra,  the notion of point value is   undertaken  in \cite{Ni}. As for the latter, as already explained, it will depend whether  this algebra is or not complete.  If it is  not complete, then  what  comes next  might be the Umweg, the roundabout route.      The original  full Colombeau algebra might be all that is needed to have a Generalized Differential  Geometry  able to handle differential equations on manifolds  having  products of distributions in their data. This is settled in   \cite{paulo}.To get  the hang of it,  we first  look at the case of the simplified algebra.


\section{Generalized Manifolds}

In this section,  $\alpha$ will stand for an index and not for our standard gauge defined in the previous section.   In \cite{OJR} the definition of a generalized manifold was given and  proved  that each such manifold had a maximal   ${\cal{G}}-$atlas. These manifolds were denoted short by ${\cal{G}}-$manifolds and the atlas by ${\cal{G}}-$atlas. In case the underlying field is $\mathbb{R}$,  we have a real generalized manifold and in case the underlying field is $\mathbb{C}$ we have a complex generalized manifold.  The topology in $\  \overline{\mathbb{R}}^n$ is the sharp topology and differentiability is in the sense of   Generalized Differential Calculus (\cite{OJRO, OJR}). Other notions   such as diffeomorphism and  continuity are those defined in \cite{OJRO, top1, top2, OJR}. 

One could ask the reason for pursuing   this path. First of all,  our aim is to establish the foundations of   a   Generalized Differential  Geometry based on the   Generalized Differential Calculus  developed in \cite{OJRO, OJR}. The latter has shown to be a natural extension of Newton's and Schwartz's  Calculus. Second, although there exists  an invariant Colombeau algebra (see \cite{MK3}), it is not clear yet if this is a topological algebra and, if it were, how its topology  interacts with the topology of a classical manifold. It is also not clear if this algebra is complete since continuity of the parameters involved in its construction is imposed (necessary).   Our proposal is that having the Colombeau full algebra at hand,  might  be what is needed  to  settle questions and problems on manifolds whose data involves products  of distributions.  For the sake of completeness, we recall the definition of a generalized manifold  in case the underlying field is $\ \mathbb{R}$.

\begin{defini} 

Let $M$  be a non-void  set.  A $\mathcal{G}$-atlas of dimension $n$ and class $C^{\infty}$  of $\ M$   is a family  $\mathcal{A}=\left\{(\mathcal{U}_{\lambda},\varphi_{\lambda})\right\}_{\lambda\in\Lambda}$  verifying the following conditions:

\begin{enumerate}
\item For every $\lambda\in \Lambda$ the map $\varphi_{\lambda}:\mathcal{U}_{\lambda} \longrightarrow  \overline{\mathbb{R}}^{n}$ is a bijection of the open subset $\emptyset \neq \mathcal{U}_{\lambda}\subset M$ onto  the open subset $\varphi_{\lambda}(\mathcal{U}_{\lambda}) \subset \overline{\mathbb{R}}^{n}$.

\item   $M=\displaystyle\bigcup_{\lambda\in \Lambda}\mathcal{U}_{\lambda}$

\item For every pair $\alpha, \beta \in \Lambda$, with  $\mathcal{U}_{ \alpha, \beta} = \mathcal{U}_{\alpha}\cap \mathcal{U}_{\beta} \neq \emptyset$,  the subsets $\varphi_{\alpha}(\mathcal{U}_{ \alpha, \beta})$ and $\varphi_{\beta}(\mathcal{U}_{ \alpha, \beta})$ are open contained in $\overline{\mathbb{R}}^{n}$ such that $\varphi_{\beta} \circ \varphi_{\alpha}^{-1}: \varphi_{\alpha}(\mathcal{U}_{ \alpha \beta}) \longrightarrow \varphi_{\beta}(\mathcal{U}_{ \alpha, \beta})$  is a diffeomorphism of class $\ C^{\infty}$.
\end{enumerate}

\begin{itemize}

\item  The pair $(\mathcal{U}_{\lambda},\varphi_{\lambda})$ is denominated  a local chart (or  coordinate system) of  $\mathcal{A}$.

\item  If  $\ \mathcal{U}\subset M$  and  $\varphi: \mathcal{U} \longrightarrow u(\mathcal{U}) $ is a  homeomorphism of  $\ \mathcal{U}$, where  $\varphi(\mathcal{U}) $ is an open set of $\overline{\mathbb{R}}^{n}$, the pair $(\mathcal{U},\varphi)$ is said to be compatible with $\mathcal{A}$ if for each pair $(\mathcal{U}_{\lambda},\varphi_{\lambda}) \in \mathcal{A}$ with $\mathcal{W}_{\lambda}=\mathcal{U}\cap \mathcal{U}_{\lambda}\neq \emptyset$ we have  that  $\varphi \circ \varphi_{\lambda}^{-1}: \varphi_{\lambda}(\mathcal{W}_{\lambda}) \longrightarrow \varphi(\mathcal{W}_{\lambda})$ is a   diffeomorphism of class $C^{\infty}$,  where $\varphi_{\lambda}(\mathcal{W}_{\lambda})$  and  $\varphi(\mathcal{W}_{\lambda}) $ are open subsets of  $\overline{\mathbb{R}}^{n}$.
\end{itemize}

By Zorn's Lemma, it  follows that, for a given $\mathcal{G}$-atlas $\mathcal{A}$ of dimension  $n$ on $M$, there exists an unique  maximal $\mathcal{G}$-atlas  $\mathcal{A}^{*}$ of dimension $n$ on $M$ defined by:

\begin{center}
$(\mathcal{U},\varphi)\in \mathcal{A}^{*}$, if and only if, $(\mathcal{U},\varphi) \in \mathcal{A}$  \ or \ $(\mathcal{U},\varphi)$ is compatible with $\mathcal{A}.$
\end{center}

\end{defini}

A   Generalized Manifold, or ${\cal{G}}-$manifold,  is a set $M$ with a ${\cal{G}}-$atlas defined on it.   A ${\cal{G}}-$manifold $M$ with a maximal ${\cal{G}}-$atlas is called a ${\cal{G}}-$differential structure of $\ M$.  If  clear from the context, the prefix ${\cal{G}}$ will be omitted. The topology on the  ${\cal{G}}-$manifold is the one that makes all charts simultaneously homeomorphisms. Our first step in setting the foundations of  the  Generalized Differential Geometry is   settling  the invariance of the dimension of a ${\cal{G}}-$manifold. 

\begin{teo}{[Dimension Invariance \cite{jose}]}
Let $(M,\mathcal{A})$ be a $\mathcal{G}$-manifold. Then, the dimension of  a $\mathcal{G}$-atlas $\mathcal{A}$  is constant in each connected component of $M$.

\end{teo}

\begin{proof}  Suppose there are two intersecting  local  charts  $(\mathcal{U}_{\alpha}, \phi_{\alpha})$ and $(\mathcal{U}_{\beta}, \phi_{\beta})$ belonging to the   $\mathcal{G}$-atlas $\mathcal{A}=\left\{(\mathcal{U}_{\alpha}, \phi_{\alpha})\right\}_{\alpha\in\Lambda}$,  such that   $\phi_{\alpha}(\mathcal{U}_{\alpha})  \subset  \overline{\mathbb{R}}^{m}$,  and  $\phi_{\beta}(\mathcal{U}_{\beta})  \subset  \overline{\mathbb{R}}^{n}$.   If  $\ \mathcal{U}_{\alpha, \beta }=\mathcal{U}_{\alpha} \cap \mathcal{U}_{\beta}\neq \emptyset$, then we have  that   $\phi_{\beta}\circ \phi^{-1}_{\alpha}: \phi_{\alpha}(\mathcal{U}_{\alpha \beta }) \longrightarrow \phi_{\beta}(\mathcal{U}_{\alpha \beta })$  is a diffeomorphism  and therefore its differential in a point   $ p \in \phi_{\alpha}(\mathcal{U}_{\alpha \beta })$, $D(\phi_{\beta}\circ \phi^{-1}_{\alpha})(p): \overline{\mathbb{R}}^{m} \longrightarrow \overline{\mathbb{R}}^{n}$,  is a $\overline{\mathbb{R}}-$isomorphism of  $\overline{\mathbb{R}}-$modules. Since $\overline{\mathbb{R}}$ is a commutative ring with unity,  it follows from  a result of  \cite{atiyah} that    $n=m$. 
\end{proof}

 Of  course,  we can define such manifolds for any differential class $C^k,  k\geq 2$. The definition is  the same, only  substituting $\ C^{\infty}$ by $\ C^{k}$. Only one  explicit, but rather trivial, example of a ${\cal{G}}-$manifold was given in \cite{OJR}.  We shall prove that any  classical  $n-$dimensional manifold $M$ can be  discretely embedded in a ${\cal{G}}-$manifold  of the same dimension thus proving that Generalized Differential  Geometry is a natural extension of Classical Differential  Geometry with the advantage that problems involving  irregular data can be considered translating them into the ${\cal{G}}-$context.

 Let $\cal{A}=\{ (U_{\alpha},\phi_{\alpha}), \alpha\!\in$$\Lambda\}$,  be an altlas of a  $C^{\infty}$  $n$-dimensional connected  submanifold $M$ of $\ \mathbb{R}^N$.   Suppose that for each $\alpha\in \Lambda$ we have that $\phi_{\alpha}(U_{\alpha})=\Omega_0 = B_r (0)\subset \mathbb{R}^n$,  the open ball of fixed radius $r >0$  centered at the origin.
 
 Denote by $\widetilde{M}_c$  the subset of  $\ \widetilde{\mathbb{R}}^N_c\subset \overline{\mathbb{R}}^N$ constructed from    $\ M$ (see the previous sections). We  saw that  $M$ is discretely  embedded in   $\widetilde{M}_c$ as constants nets.   Recall from \cite{OJR} that  $\widetilde{\mathbb{R}}^N_c \subset \overline{B}_1 (0)$, the  ball of  radius $ 1$  centered at the origin and that the image of $\mathbb{R}^N$ under this map is a grid of  equidistant  points.  We denote by  $\widetilde{\Lambda}$  the set of   maps from  $\ I =]0,1]$ into $\ \Lambda$ and, for $\lambda \in \widetilde{\Lambda}$,  we  denote by  
 
 $$U_{\lambda}=\  \mbox{ the strongly internal set  \ }\langle U_{\lambda(\varepsilon)}\rangle\  \mbox{contained in}\   \widetilde{\mathbb{R}}^N_c $$ 
 
 $$  \phi_{\lambda}=(\phi_{\lambda(\varepsilon)})_{\varepsilon \in I}$$
 
 $$  \phi_{\lambda}\ : U_{\lambda}\longrightarrow \widetilde{\mathbb{R}}^n_c, \mbox{defined by} \   \phi_{\lambda}([(p_{\varepsilon})])=[(\phi_{\lambda(\varepsilon)}(p_{\varepsilon}))] $$

For  $p=[(p_{\varepsilon})]\in\widetilde{M}_c$,  consider   the set $\{q\in M\ : \exists\  \varepsilon_n\rightarrow 0,  \ p_{\varepsilon_n}\rightarrow q\}$.  Algebraically this can be written as:  Given  $q_0\in \mathbb{R}^n$, we have that   $\ q_0\in \{q\in M\ : \exists\  \varepsilon_n\rightarrow 0,  \ p_{\varepsilon_n}\rightarrow q\}$ if and only if there exists $e\in {\cal{B}}(\overline{\mathbb{R}})$ such that $e\cdot p\approx e\cdot q_0$ (extending the notion of association to $\overline{\mathbb{K}}^n$ in the obvious way). This is a compact subset of $\ M$ to which  we  shall refer  as the {\it{support of the point}} $p$ and denote  it  by $supp(p)$. It follows that there exists a complete set of orthogonal  idempotents  $(e_{\lambda})$ such that 

$$p\ =\sum\limits_{x_{\lambda}\in supp(p)}e_{\lambda}\cdot p  $$ $$ (e_{\lambda}\cdot p \approx e_{\lambda}\cdot x_{\lambda})$$

\noindent For example, if   $p=[(p_{\varepsilon})]=[(sin(\frac{1}{\varepsilon}))]$ then $supp(x)=[-1,1]$. On the other hand, if $p\in B_1(0)$ then $supp(p)=\{0\}$.  Although $supp(p)$ might be uncountable, nevertheless the sum above is well defined and can be thought of as   a history of events: multiplying by idempotents one sees its behavior along a specific path. It generalizes the concept of interleaving.     The support of elements   belonging to a halo of a point in  $\ \mathbb{R}^n$  consists of only that  single point.  

Our next result  gives conditions under which the objects defined above turn  $\widetilde{M}_c$ into a ${\cal{G}}-$manifold of dimension $n$.  This will enable us to look at  Generalized  Differential Geometry as a natural  extension of  Classical Differential Geometry. The latter having as basis Newton's Differential Calculus and the former having as basis the  Generalized Differential Calculus defined and developed in \cite{OJRO, OJR}.

\begin{prop} 
\label{condicao}
Suppose that for each $\alpha\in \Lambda$ the map $\alpha$ and its inverse are Lipschitz functions with respect to the norms of $\ \mathbb{R}^N$ and $\ \mathbb{R}^n$. Then the following hold.

\begin{enumerate}

\item The topology of $\  \widetilde{M}_c$ is induced by the topology of $\ \overline{\mathbb{R}}^n$.
\item If  $\lambda\in \widetilde{\Lambda}$ has finite range then $\ U_{\lambda}$ is an open subset of $\ \widetilde{M}_c$.
\item $ \langle \Omega_0\rangle$ is an open subset of $\  \widetilde{\mathbb{R}}^n_c$
\item If  $\lambda\in \widetilde{\Lambda}$ has finite range then   $  \phi_{\lambda}\ : U_{\lambda}\longrightarrow\langle \Omega_0\rangle\subset \widetilde{\mathbb{R}}^n$ is a isometry w.r.t  to the sharp topologies  of  $\ \overline{\mathbb{R}}^n$ and $\ \overline{\mathbb{R}}^N$.
\item  $\widetilde{M}_c= \bigcup\limits_{\lambda \in \widetilde{\Lambda}}U_{\lambda} $, with $\lambda$ of finite range.
\item If  $\lambda_1 ,\lambda_2\in \widetilde{\Lambda}$ are of finite range and $U_{\lambda_1,\lambda_2}:=U_{\lambda_1}\cap U_{\lambda_2}\neq \emptyset$ then $\phi_{\lambda_2}\circ \phi_{\lambda_1}^{-1}$ is a $C^{\infty}$ diffeomorphism on its domain $\phi_{\lambda_1}(U_{\lambda_1,\lambda_2})$.

\end{enumerate}

\end{prop}

\begin{proof} Let   $p=[(p_{\varepsilon})]\in\widetilde{M}_c$. Since $supp(p)$ is compact,    there exist  a finite number of    $\ \alpha\in \Lambda$ such that   $supp(p)$ is contained in the union of the corresponding $\ U_{\alpha}$'s.  Let $\delta$ be a Lebesgue number of this covering and choose a finite number of $q_i\in supp(p)$ such that  $supp(p) \subset \bigcup\limits_{i} B_{\delta_1}(q_i)$,  with $\delta_1=\frac{\delta}{2}$. Starting with $q_1$, define $\lambda(\varepsilon)=\alpha_{q_1}\in \Lambda$, where $\alpha_{q_1}$ is such that $B_{\delta_1}(q_1)\subset U_{\alpha_{q_1}}$  and $p_{\varepsilon}\in B_{\delta_1}(q_1)$.

  For $\lambda(\varepsilon)$ not yet defined,  continue defining $\lambda(\varepsilon)= \alpha_{q_2}\in \Lambda$,  where $\alpha_{q_2}$ is such that $B_{\delta_1}(q_2)\subset U_{\alpha_{q_2}}$ and $p_{\varepsilon}\in B_{\delta_1}(q_2)$. Since there are a finite number of $q_i$'s, this process ends in a finite number of steps. If $\lambda$ is not defined on $I$ then there exists a sequence $\  (\varepsilon_n)$,  converging to $0$, with  $p_{\varepsilon_n}\rightarrow q\in supp(p)$. Since the balls $B_{\delta_1}(q_i)$ cover $A_p$, there exists a $n_0$ such that $n>n_0$ implies that $p_{\varepsilon_n}$ is in the ball $\ B_{\delta}(q_{i_0})$, say. But this is a contradiction, since,  for these $\varepsilon$'s, $\lambda(\varepsilon)$ was already defined. Hence, we defined a $\ \lambda\in \widetilde{\Lambda}$ such that $p_{\varepsilon}\in U_{\lambda(\varepsilon)}, \varepsilon \in I$ and the distance to the boundary is bigger than $\frac{\delta}{2}$.  It follows that $p\in U_{\lambda}$, with $\ \lambda$  being of finite range.

 For each $\lambda(\varepsilon)$ there exists an open subset $U^{\lambda(\varepsilon)}\subset \mathbb{R}^N$ such that $U_{\lambda(\varepsilon)}=M\cap U^{\lambda(\varepsilon)}$.  Setting  $U^{\lambda}= \langle U^{\lambda(\varepsilon)}\rangle$, we have that  $U_{\lambda}=\widetilde{M}_c \cap U^{\lambda}$, with $U^{\lambda}$ an open subset of $\ \overline{\mathbb{R}}^N$, proving that $\widetilde{M}_c$ has the induced topology.  The fact that  strongly internal sets are open can be found in  \cite{OV}.  This  settles the proof of the  first three and the fifth  items.

To finish the proof, we prove the forth and sixth items. We first prove that $\phi_{\lambda}$ is well defined. In fact, since local charts are Lipschitz,  and  the $\lambda$'s are  of finite range, it follows easily  that  $\|p-q\| = \|\phi_{\lambda} (p)-\phi_{\lambda}(q)\|$, i.e., the $\phi_{\alpha}$'s are isometries considered as maps from a subset of $\widetilde{\mathbb{R}}^N$ to a subset of  $\ \widetilde{\mathbb{R}}^n$.  From this it follows that if   $dist(p, (U_{\lambda})^c)$ is an invertible then $dist(\phi_{\alpha}(p), (\langle \Omega_0\rangle)^c)$ is an invertible, thus proving that $\phi_{\lambda}$ is well defined and its image is in $U_{\lambda}$. Surjectivity,  injectivity and continuity are now obvious.   The  map resulting from a  change of coordinates,  is a homeomorphism that  stems from a net whose elements are all infinitely differentiable taking value in a bounded subset of $\ \mathbb{R}^n$.  Hence it originates a diffeomorphism.    \end{proof}

The condition that all charts have the same image  is not necessary. We could just suppose that there exists $r>0$ such that $B_r(0)\subset \ \mathbb{R}^n$ contains all of them. The $(U_{\lambda},\phi_{\lambda})$ are called {\it{local charts }} of $\ \widetilde{M}_c$. If $\lambda$ is constant then  $U_{\lambda}$ is called a {\it{principal}} chart. For  $\lambda\in \widetilde{\Lambda}$ and  $\  x\in B_1(0)$ define $x\lambda$ as  $x\lambda(\varepsilon)=\lambda(x(\varepsilon))$. For an idempotent $e\in {\cal{B}}(\overline{\mathbb{K}})$ we define $e\lambda=e(\varepsilon)\lambda(\varepsilon)$, meaning that when $e(\varepsilon)=1$ this index will be omitted. With this notation fixed, the proof of the proposition gives  us the following corollary.

\begin{cor} Let $U_{\lambda}$ be a local chart with $\ \lambda$ of finite range. There exist principal charts $U_{\lambda_i}, i= 1,k$  and a complete set of mutually orthogonal idempotents $e_1,\cdots, e_k$ such that  $U_{\lambda}\subset \bigcup U_{\lambda_i}$ and $U_{\lambda}= \bigcup U_{e_i\lambda_i}$ 
\end{cor}

The corollary shows  the importance of principal charts since they can serve as tools to use sheaf like arguments in the proofs.  For results to hold, they just have to be proved for principle charts. 

\begin{teo}  Let $M$ be a submanifold of $\ \mathbb{R}^N$ of  dimension $n$, and suppose that  its local charts and their   inverses are Lipschitz functions with respect to the norms of $\ \mathbb{R}^N$ and $\mathbb{R}^n$. Then  $(\widetilde{M}_c,{\cal{A}})$, ${\cal{A}}=\{ (U_{\lambda},\phi_{\lambda})\ :\ \lambda \in \widetilde{\Lambda} \ \mbox{of finite range}\}$ is a generalized submanifold  of $\ \overline{\mathbb{R}}^N$ of codimension $N-n$ containing $\ M$ as a discrete subset. Moreover, each local chart is an isometry in the sharp topologies and the geometry of $\widetilde{M}_c$ extends in a natural way the geometry  of $\ M$.
\end{teo}

\begin{proof}
Only the last part of the theorem needs to be proved. To see this,  we  use Lemma A.1 of  \cite{V2} which state that for  each compact subset of $K\subset \ M$ its  Riemannian metric  satisfies a $\|p-q\|\leq dist_M(p,q)\leq C\cdot \|p-q\|$, for some $C>0$, $p,q\in K$ and $dist_M$ the Riemannian metric of $\ M$. This implies that local charts are isometries.\end{proof}

Using Whitney's  Embedding Theorem,  the above theorem extends to abstract manifolds.     The construction of an atlas for $\widetilde{M}_c$ shows that if $\ M$ is a classical manifold   then  $\widetilde{M}_c$ can be covered with a countable number of local charts. Another important fact in the construction of  $ \widetilde{\Lambda}$ is that since its elements have finite range and we can work with  local charts $(U, \phi)$ of $\ M$  such that $U$ is relatively compact in $\ M$, it follows that for $\lambda_1,\lambda_2\in  \widetilde{\Lambda}$ we have that $\det(D(\phi_{\lambda_2}\circ\phi_{\lambda_1}^{-1})(p))\in Inv(\overline{\mathbb{R}})$, for each point $\ p$ in its domain. If we define orientabillity   as is done the classically, this shows that if $\ M$ is orientable then so is $\ \widetilde{M}_c$. In particular, complex  ${\cal{G}}-$manifolds will be orientable. It is   important to note that the construction of the differential structure does not depends on the Riemannian metric of the manifold $\ M$ (see \cite{MK4}).  

Using the Generalized  Fixed Point Theorem, one can prove the global existence of geodesics without  using the usual stratagem: proving existence depending on $\ \varepsilon$ and then  proving moderateness.  There is no need to expand on this since everything works just as in the classical case.  The environments  of Classical Differential Geometry   are    discretely embedded in the the environments  of Generalized Differential Geometry.


\section{Differential Functions on ${\cal{G}}-$manifold}

A one dimensional ${\cal{G}}-$manifold   shall be referred  to as  a {\it curve}. This definition agrees with the notion of a history given in \cite{OJR}, only now we parametrize  the history.  If $M$ is an $n-$dimensional  ${\cal{G}}-$manifold  and $p\in M$ then the tangent vectors  and  space at $p$  are  defined just as in the classical way (see \cite{jose}). The tangent space we shall also   denote by $T_pM$.  Furthermore, it is easily proved that $T_pM$  is a free  $\overline{\mathbb{R}}-$module which is $\overline{\mathbb{R}}-$isomorphic to $\overline{\mathbb{R}}^n$.

As in the classical case,  one  defines  the tangent bundle of $\ M$ and denote it by $\ TM$. It is easily seen that if $\ M$ is a classical $n$-dimensional manifold then $\widetilde{(TM)}_c=T\widetilde{M}_c$ is a $\ {\cal{G}}-$manifold of dimension $n^2$. 

Given another  ${\cal{G}}-$manifold  $\ N$,    define a differential function between $\ M$ and $\ N$ using the same classical definition.  In case $\ N=\overline{\mathbb{R}}$ we call a differentiable function a scaler field  and if  $\ N=\overline{\mathbb{R}}^n$, with $n>1$,   then  we call such a map a vector valued map, being a vector field if also the dimension of $\ M$ equals $n$.    

The notions of an  immersion and an  embedding are defined completely analogous as in classical geometry.  Using the Chain rule of Colombeau Generalized Calculus it follows easily that composition of  differentiable maps between  ${\cal{G}}-$manifolds also satisfy the Chain Rule (see \cite{jose}). We recall a   Linear Algebra results  (see   \cite{grosser1}). 

\begin{lema}

Let  $A:\overline{\mathbb{R}}^{n}\longrightarrow \overline{\mathbb{R}}^{n}$ be a $\overline{\mathbb{R}}-$linear map. Then $A$ is injective if and only if it is surjective if and only if $\det(A)\in Inv(\overline{\mathbb{R}})$.
\end{lema}

We sum up, without proofs, some of the most classical theorems that also hold for ${\cal{G}}-$manifolds. Some of the  proofs rely on the previous lemma and the fact  that $Inv(\overline{\mathbb{R}})$ is  open (see \cite{jose}).

\begin{teo}

Let $M_{1},M_{2}$ and $M_{3}$ be $\mathcal{G}$-manifolds. If $f:M_{1} \longrightarrow  M_{2}$ and $g:M_{2} \longrightarrow M_{3}$ are differentiable applications at $p\in M_{1}$ and $f(p)\in M_{2}$, respectively, then  $g \circ f:M_{1} \longrightarrow M_{3}$ is differentiable at $p$ and $D(g \circ f)_{p}=(Dg)_{f(p)}\circ Df_{p}.$
\end{teo}

\begin{teo} 

Let $M_{1}$ and $M_{2}$ be  $n$-dimensional  $\mathcal{G}$-manifolds and  $f: M_{1} \longrightarrow M_{2}$ a map   of  class $C^{\infty}$, such that for   $p_0 \in M_{1}$ we have that $Df_{p_0}:T_{p_0}M_{1}\longrightarrow T_{f(p_0)}M_{2}$, is an isomorphism.  Then  $f$ is a local diffeomorphism of class $C^{\infty}$.

\end{teo}

\begin{teo} 

Let $f: M\longrightarrow N$ be  an immersion at $p$ of class $C^{\infty}$, where $M$ and $N$ are  $\mathcal{G}$-manifolds of dimension $m$ and $n$, respectively.   Then   there exist   local coordinate systems around $p$ and $f(p)$, such that
\begin{center}
$f(x_{1},...,x_{m})=(x_{1},...,x_{m},0,...,0).$
\end{center}
\end{teo}

\begin{teo}

If $f: M \longrightarrow N$ is an generalized embedding, then $f(M)$ is an $\mathcal{G}$-submanifold of $N$.
\end{teo}

\begin{defini}
Let $f: \widetilde{\Omega}_{c}  \subset \overline{\mathbb{R}}^{n} \longrightarrow  \overline{\mathbb{R}}^{m}$ be a differentiable map,  where $\Omega$ is an open subset  of  $\mathbb{R}^{n}$.  A point $a \in \overline{\mathbb{R}}^{m}$ is called a  regular value of $f$  if  for each $x \in f^{-1}(a)$ the  derivative $f'(x): \overline{\mathbb{R}}^{n} \rightarrow \overline{\mathbb{R}}^{m}$  is   surjective. 
\end{defini}

\begin{teo} 
Let $\Omega$ be an open subset of $\mathbb{R}^{m} \times \mathbb{R}^{n}$ and  $f: \widetilde{\Omega}_{c}  \longrightarrow  \overline{\mathbb{R}}^{n}$ be a  application of class $\mathcal{C}^{\infty}$, where $\widetilde{\Omega}_{c} \subset  \overline{\mathbb{R}}^{m} \times  \overline{\mathbb{R}}^{n}$.   If  $a \in Im(f)$ is a regular value of $f$, then:

\begin{enumerate}

\item   $f^{-1}(a) $ is an $m-$dimensional  $\mathcal{G}$-submanifold of  $\ \overline{\mathbb{R}}^{m} \times  \overline{\mathbb{R}}^{n}$.
\item For each $p \in f^{-1}(a)$, we have that $T_{p}(f^{-1}(a)) = \ker( f^{'}(p))$.
 \end{enumerate}
\end{teo}

Let's look at some examples of $\  {\cal{G}}-$manifolds.

\begin{ex}
$\ $\newline

\begin{enumerate}

\item Consider $M= Graf (f)$, where $f \ : \Omega\subset \mathbb{R}^n \longrightarrow \mathbb{R}$, a $C^{\infty}-$function with bounded first derivate. Denote by $\phi : M\longrightarrow \Omega$ the projection $\phi(p)=x$, where $p=(x,f(x)) \in M$. If $q=(y,f(y))\in M$ then $\| \phi(p)-\phi(q)\|=\|x-y\|<\|p-q\|$. On the other hand   

$\|\phi^{-1}(x)-\phi^{-1}(y)\|^2  =\|(x,f(x))-(y,f(y))\|^2=\|x-y\|^2+|f(x)-f(y)|^2\leq \|x-y\|^2+\|\nabla f(p_0)\|\cdot \|x-y\|)^2\leq (1+C)\cdot \|x-y\|^2$

This proves that the conditions of Proposition~\ref{condicao}  are satisfied and thus $\widetilde{M}_c$ is a ${\cal{G}}-$submanifold of $\overline{\mathbb{R}}^{n+1}$.

\item  Let $M\subset \mathbb{R}^n$ be a codimension one submanifold   with an atlas whose elements are graphs. Then, by Theorem 3 of the previous section,   we have that  $\widetilde{M}_c$ is a ${\cal{G}}-$submanifold of $\ \overline{\mathbb{R}}^{n+1}$.  Hence, this is true if $M$ is a $m-$dimensional surface of $\ \mathbb{R}^n$. In particular, this holds if $M$ is the pre image of  a regular value of a $C^{\infty}$ differentiable function $f :\mathbb{R}^n\longrightarrow \mathbb{R}$.

\item Let  $M=S^n_r \subset \mathbb{R}^n$ be an  $n-$dimensional sphere of radius $r$.  It can be parametrized by graphs such that derivates of the functions involved are bounded. Hence, by the previous example, we have that $\widetilde{M}_c$ is a ${\cal{G}}-$submanifold of $\ \overline{\mathbb{R}}^{n+1}$. One can cover $\ M$ just with two local charts but this can not be done with $\ \widetilde{M}_c$. This example inspired the construction of the non principal  charts and the notion of the support of a generalized point seen in  the previous section.

\item The sphere $S=S_1(0)$ contained in $\ \overline{\mathbb{R}}^n$ is a generalized manifold whose local charts do not come from subsets of $\ \mathbb{R}^n$. In fact,  $S$ is an open subset of $\ \overline{\mathbb{R}}^n$ because given $x\in$ and $y\in B_1(0)$ we have that $\|x+y\|=max\{\|x\|, \|y\|\}=1$, since $\|x\|=1>\|y\|$. Consequently, we can take local charts to be the identity map with domain $B_1(x)$. Since these balls are either equal or disjoint, it follows that they form an atlas for $S$.

\item  Let $f\in {\cal{G}}(\Omega)$ with $\Omega\subset \mathbb{R}^n$. We know (see \cite{OJR})  that $f$ can be viewed as a differentiable map $\widetilde{\Omega}_c\longrightarrow \overline{\mathbb{R}}$ and its differential at each point is a $ \overline{\mathbb{R}}-$linear map from  $\overline{\mathbb{R}}^n$ to $\overline{\mathbb{R}}$. A value $a\in Im(f)$ is said to be {\it{a regular value}} of $f$ if for each $x\in f^{-1}(a)$ we have that $Df$ is surjective. This only happens if,  writing  $\nabla f(x)=(z_1,\cdots, z_n)$,  the ideal in $\  \overline{\mathbb{R}}$ generated by $z_1,\cdots, z_n$ equals $\overline{\mathbb{R}}$. In particular,  this is the case if $\ \|\nabla f(x)\|_2^2  \in Inv(\overline{\mathbb{R}})$. 
If $a$ is a regular value of $\ f$ set $M=f^{-1}(a)$. We assert that $M$ is a submanifold of $\ \overline{\mathbb{R}}^n$. In fact, just like in the classical case, we can use the Implicit Function Theorem (see \cite{OJRO}) to prove that at every point $M$ is locally a graph over a subset of $\widetilde{\Omega}_c$.  The standard classical argument still holds to complete the assertion.

An example of such a function is $\ f(x)=\|x\|_2^2$ and the value in question is $a=1$. In this case we have that $\|\nabla f(x)\|^2_2=4\|x\|_2^2=4\in Inv(\overline{\mathbb{R}})$. 

\item The  halo of  any point in $\ \overline{\mathbb{K}}^n$  is a generalized manifold which is not classical.  In \cite{OJR} an example of a function $\ f\not \equiv 0$ is given such that $f^{\prime}\equiv 0$. This function is $\ f(x)=\alpha_{-2\ln(\|x\|)}$ which is not a Colombeau  generalized function.  This $\ f$ is constant on spheres $\ S_R(0)$ and the only point where it is not locally constant is $\ x_0=\vec{0}$. This  $\ f$ is easily modified  such that it is of class $\ C^k$. Since the origin is the only points where such spheres accumulate, and spheres are clopen,  we have that $\ Graf(f)-\{\vec{0}\}$ is a generalized manifold but $\ Graf(f)$ is not a generalized manifold. 
\end{enumerate}

\end{ex}


  \section*{Acknowledgements}
  
 The second author is indebted to the University of S\~{a}o Paulo for its warm hospitality and to   the  Federal University of Roraima,  for providing  the opportunity to do this research.


\Addresses

\end{document}